\documentclass[11pt]{amsart}

\usepackage[utf8]{inputenc}
\usepackage[a4paper,twoside]{geometry}

\usepackage{enumitem}
\usepackage{amsmath, amssymb, amsthm}

\usepackage[english]{babel}
\usepackage{fancyhdr}

\usepackage[colorlinks=true,linkcolor=blue,citecolor=red,backref=page,pdftex]{hyperref}

\newcommand{\retrait}{\hspace{1.7cm}}

\newcommand{\unp}{\mathbf{\mathrm{1 \kern-0.25em I}}}
\newcommand{\un}{\mathbf{1}}

\newcommand{\X}{\mathbf X}
\newcommand{\Hb}{\mathbf H}

\newcommand{\V}{\mathbb V}
\newcommand{\W}{\mathbb W}
\newcommand{\R}{\mathbb R}
\newcommand{\N}{\mathbb N}

\newcommand{\Q}{\mathbb Q}
\newcommand{\G}{\mathbf G}
\newcommand{\Z}{\mathbb Z}

\newcommand{\T}{\mathbb T}

\newcommand{\di}{{\rm d}}

\newcommand{\cal}{\mathcal}

\newcommand{\esp}{\mathbb E}
\newcommand{\prob}{\mathbb P}

\newcommand{\supp}{\mathrm{supp}\,}

\newtheorem{theorem}{Theorem}[section]
\newtheorem*{theorem*}{Theorem}
\newtheorem*{lemma*}{Lemma}
\newtheorem*{corollary*}{Corollary}
\newtheorem{lemma}[theorem]{Lemma}
\newtheorem{proposition}[theorem]{Proposition}
\newtheorem*{proposition*}{Proposition}
\newtheorem{corollary}[theorem]{Corollary}

\theoremstyle{definition}
\newtheorem{definition}[theorem]{Definition}
\newtheorem*{definition*}{Definition}
\newtheorem{example}[theorem]{Example}

\newtheorem*{notations*}{Notations}
\theoremstyle{remark}
\newtheorem{remark}[theorem]{Remark}

\numberwithin{equation}{section}
\begin{document}

\title{On the affine random walk on the torus}


\author{Jean-Baptiste Boyer}
\email{maths@jbboyer.fr}


\date{\today}

\begin{abstract}
Let $\mu$ be a borelian probability measure on $\mathbf{G}:=\mathrm{SL}_d(\mathbb{Z}) \ltimes \mathbb{T}^d$. Define, for $x\in \mathbb{T}^d$, a random walk starting at $x$ denoting for $n\in \mathbb{N}$,
\[
\left\{\begin{array}{rcl}
X_0 &=&x\\
X_{n+1} &=& a_{n+1} X_n + b_{n+1}
\end{array}\right.
\]
where $((a_n,b_n))\in \mathbf{G}^\mathbb{N}$ is an iid sequence of law $\mu$.

Then, we denote by $\mathbb{P}_x$ the measure on $(\mathbb{T}^d)^\mathbb{N}$ that is the image of $\mu^{\otimes \mathbb{N}}$ by the map $\left((g_n) \mapsto (x,g_1 x, g_2 g_1 x, \dots , g_n \dots g_1 x, \dots)\right)$ and for any $\varphi \in \mathrm{L}^1((\mathbb{T}^d)^\mathbb{N}, \mathbb{P}_x)$, we set $\mathbb{E}_x \varphi((X_n)) = \int \varphi((X_n)) \mathrm{d}\mathbb{P}_x((X_n))$.

Bourgain, Furmann, Lindenstrauss and Mozes studied this random walk when $\mu$ is concentrated on $\mathrm{SL}_d(\mathbb{Z}) \ltimes\{0\}$ and this allowed us to study, for any hölder-continuous function $f$ on the torus, the sequence $(f(X_n))$ when $x$ is not too well approximable by rational points.

In this article, we are interested in the case where $\mu$ is not concentrated on $\mathrm{SL}_d(\mathbb{Z}) \ltimes \mathbb{Q}^d/\mathbb{Z}^d$ and we prove that, under assumptions on the group spanned by the support of $\mu$, the Lebesgue's measure $\nu$ on the torus is the only stationary probability measure and that for any hölder-continuous function $f$ on the torus, $\mathbb{E}_x f(X_n)$ converges exponentially fast to $\int f\mathrm{d}\nu$.

Then, we use this to prove the law of large numbers, a non-concentration inequality, the functional central limit theorem and it's almost-sure version for the sequence $(f(X_n))$.

\medskip
In the appendix, we state a non-concentration inequality for products of random matrices without any irreducibility assumption.
\end{abstract}

\maketitle

\tableofcontents

\section{Introduction and main results}

Let $d\in \N$, $d\geqslant 2$ and $\T^d:= \R^d/\Z^d$ be the torus of dimension $d$. Let $\mu$ be a borelian probability measure on $\G:=\mathrm{SL}_d(\Z) \ltimes \T^d$. Define, for any $x\in \T^d$, a random walk starting at $x$ by denoting, for any $n\in \N$,
\[
\left\{\begin{array}{rcl}
X_0 &=&x\\
X_{n+1} &=& a_{n+1} X_n + b_{n+1}
\end{array}\right.
\]
where $((a_n,b_n))\in \G^\N$ is an iid sequence of law $\mu$.

Then, we denote by $\prob_x$ the measure on $\X^\N$ that is the image of the measure $\mu^{\otimes \N}$ on $\G^\N$ by the map $\left((g_n) \mapsto (x,g_1 x, g_2 g_1 x, \dots , g_n \dots g_1 x, \dots)\right)$ and by $\esp_x$ the operator of integration against the measure $\prob_x$.

We denote by $P$ the Markov operator associated to  $\mu$. This is the operator defined for any borelian non-negative function $f$ on $\T^d$ and any $x\in \T^d$ by
\[
Pf(x) = \int_\G f(gx)\di\mu(g)
\]
Thus, for any $n\in \N$, we have that
\[
P^n f(x) = \int_\G f(gx) \di\mu^{\ast n}(g) = \int_{\T^d} f(y) \di\mu^{\ast n} \ast \delta_x(y) = \int_{(\T^d)^\N} f(X_n) \di\prob_x((X_n)) = \esp_x f(X_n)
\]
Where we noted $\mu^{\ast n}$ the $n-$th power of convolution of the measure $\mu$ ($\mu^{\ast 0}$ is by convention the Dirac measure at $(I_d,0)$).

\bigskip
Bourgain, Furmann, Lindenstrauss and Mozes studied in~\cite{BFLM11} the case where $\mu$ is concentrated on $\mathrm{SL}_d(\Z) \ltimes \{0\}$ and they proved that, under assumptions on the support of $\mu$, the only $P-$invariant probability measures on the torus where the Lebesgue's measure $\nu$ and the uniform measures on unions of rational orbits (which are finite). their result is even more precise since they give the rate of convergence of $\esp_x f(X_n)$ to $\int f\di\nu$ in terms of diophantine properties of $x$ and this allowed us to study the sequence $(f(X_n))$ for starting points $x$ that are not too well approximable by rational points in~\cite{Boytore}.

\medskip
In this article, we are interested in the case where $\mu$ is not concentrated on $\mathrm{SL}_d(\Z)\ltimes \Q^d/\Z^d$. A result by Benoist-Quint (see~\cite{BQStat0}) shows that in this case, under assumptions on the projection on $\mathrm{SL}_d(\Z)$ of the subgroup spanned by the support of $\mu$, the only $P-$invariant probability measure on the torus is Lebesgue's measures and this proves that for any continuous function $f$ on $\T^d$ and any $x\in \T^d$,
\[
\frac 1 n \sum_{k=0}^{n-1} f(X_k) \xrightarrow\, \int f\di\nu \;\prob_x-\text{a.e.}
\]
The aim of this article is to precise the previous convergence by proving a Central Limit Theorem, a Law of the Iterated Logarithm, etc.

To do so, we are going to make a few assumptions on the subgroup spanned by $\{a|(a,b)\in\supp\mu\}$.

In the sequel, we will say that a closed subgroup $\Hb$ of $\mathrm{SL}_d(\R)$ is strongly irreducible if it doesn't fix any finite union of non-trivial subspaces of $\R^d$. Moreover, we will say that $\Hb$ is proximal if it contains an element $g$ such that there is $v_g^+ \in \R^d\setminus\{0\}$, $\lambda \in \R$ and a $g-$
invariant hyperplane $V_g^<$ in $\R^d$ such that $\R^d=  \R v_g^+ \oplus \V_g^<$, $gv_g^+ = \lambda v_g^+$ and the spectral radius of $g$ in $V_g^<$ is strictly smaller that $|\lambda|$. Finally, we will say that a probability measure $\mu$ on $\mathrm{SL}_d(\R)$ is strongly irreducible and proximal if the closure of the subgroup spanned by the support of $\mu$ has these properties.

These two assumptions are actually assumptions on the Zariski-closure of $\Hb$ and so, as an example, they are satisfied if $\Hb$ is Zariski-dense in $\mathrm{SL}_d(\R)$.

Finally, we will say that a measure $\mu$ on $\mathrm{SL}_d(\R)$ has an exponential moment if there is some $\varepsilon\in \R_+^\ast$ such that
\[
\int_{\mathrm{SL}_d(\R)} \|g\|^\varepsilon \di\mu(g)<+\infty
\]
We will see in the sequel that our study of the random walk on the torus requires arguments of orbit closures. This is why we give a name to the property that we will use and we will see right after examples of measures satisfying it.
\begin{definition}
Let $\mu$ be a borelian probability measure on $\G$.

We say that $\mu$ satisfies an effective shadowing lemma if for any $C',t' \in \R_+^\ast$, there are $C_1,C_2,M,t,L \in \R_+^\ast$ such that for any $x,y\in \T^d$, any $r\in \R_+^\ast$ and any $n\in \N$, with $r\leqslant C_1e^{-Ln}$, if
\[
\mu^{\ast n} \left(\left\{ g\in \G\middle| d(gx,y) \leqslant r \right\}\right) \geqslant C_2 e^{-t n}
\]
then, there are $x',y'\in \T^d$ such that $d(x,x'),d(y,y')\leqslant re^{Mn}$ and
\[
\mu^{\ast n} \left(\left\{ g\in \G\middle| gx'=y'\right\}\right)\geqslant C'e^{-t'n}
\]
\end{definition}

%
%

\begin{remark}
For a measure to satisfy this property means that if a lot of elements $g\in \supp\mu^{\ast n}$ send $x$ close to $y$, it is only because $x$ and $y$ are close to points of the same orbit. The name comes from the theory of hyperbolic diffeomorphisms since when $\mu=\delta_{g_0}$, saying that $\mu$ satisfies an effective shadowing lemma means that there is some constant $M$ such that for any large enough $K$, any $n\in \N$ and any $x,y\in \T^d$ with $d(g^n_0 x,y) \leqslant e^{-K n}$, there are $x',y' \in \T^d$ such that $d(x,x'),d(y,y')\leqslant e^{-(K-M)n}$ and $g^n_0 x' = y'$.
\end{remark}

\begin{example}
This is a technical definition but we will see in section~\ref{section:fermeture} a criterion (the proposition~\ref{proposition:critere_diophantien}) that allows to tell if a measure satisfies to an effective shadowing lemma and we will deduice examples from it.

In particular, we will see in example~\ref{exemple:element_diophantien} that if $b_0\in\T^1$ is such that there are $C,L \in \R_+^\ast$ such that for any $q\in \N^\ast$, $d(qb_0,0) \geqslant \frac C {q^L}$, then, any borelian probability measure $\mu$ on $\G$ whose projection on $\mathrm{SL}_d(\Z)$ is strongly irreducible, proximal and has an exponential moment and such that $F(\mu):=\{\text{coefficients of }b|(a,b)\in \supp\mu\} \subset \{0,b_0\}$ satisfies an effective shadowing lemma.

Moreover, in example~\ref{example:diophantien_generique} we will prove that if $a_1, \dots , a_N \in \mathrm{SL}_d(\Z)$ generate a strongly irreducible and proximal group, then for a.e. $b_1, \dots , b_N\in \T^d$, the measure $\mu = \frac 1 N \sum_{i=1}^N \delta_{(a_i,b_i)}$ satisfies an effective shadowing lemma.
\end{example}

For $\alpha\in ]0,1]$, we denote by $\mathcal{C}^{0,\alpha}(\T^d)$ the space of $\alpha-$hölder-continuous functions $f$ on $\T^d$ endowed with the norm
\[
\|f\|_\alpha := \|f\|_\infty + m_\alpha(f)
\]
where
\[
\|f\|_\infty := \sup_x|f(x)| \text{ and }m_\alpha(f) := \sup_{x\not=y} \frac{ |f(x)-f(y)|}{d(x,y)^\alpha}
\]
where $d$ is the distance induced by some norm on $\R^d$.

Moreover, for any two borelian probability measures $\vartheta_1, \vartheta_2$ on $\T^d$, we denote by $\mathcal{W}_\alpha(\vartheta_1, \vartheta_2)$ the Kantorovich-Rubinstein's distance of $\vartheta_1$ and $\vartheta_2$ and this is defined by
\[
\mathcal{W}_\alpha(\vartheta_1, \vartheta_2) := \sup_{\substack{f\in \mathcal{C}^{0,\alpha}(\T^d)\\ \|f\|_\alpha\leqslant 1}} \left|\int f\di\vartheta_1 - \int f\di\vartheta_2 \right|
\]

This will allow us to prove the

\begin{theorem}\label{theoreme_equidistribution}
Let $\mu$ be a borelian probability measure on $\G:=\mathrm{SL}_d(\Z)\ltimes \T^d$ that is not concentrated on $\mathrm{SL}_d(\Z) \ltimes \Q^d/\Z^d$ and satisfies an effective shadowing lemma. Note $\mu_0$ the projection of $\mu$ on $\mathrm{SL}_d(\Z)$ and assume that $\mu_0$ is strongly irreducible, proximal and has an exponential moment.

Denote by $P$ the Markov operator associated to $\mu$.

Then, the Lebesgue's measure $\nu$ on $\T^d$ is the only $P-$invariant borelian probability measure on the torus. Moreover, for any $\alpha \in \left]0,1\right]$, there are $C,t\in \R_+^\ast$ such that for any $f \in \cal C^{0,\alpha}(\T^d)$ and any $n\in \N$,
\[
\sup_{x\in \T^d}\mathcal{W}_\alpha(\mu^{\ast n} \ast \delta_x, \nu) \leqslant Ce^{-tn} \|f\|_\alpha
\]
In particular, for any $\alpha-$hölder-continuous function $f$ on the torus, there is a continuous function $g$ such that
\[
f-\int f \di\nu = g-Pg \text{ and }\|g\|_\infty \leqslant C \|f\|_\alpha
\]
\end{theorem}

\begin{remark}
We don't know if the function $g$ that we construct in this theorem is hölder-continuous.
\end{remark}

This theorem will allow us to prove a few of the classical results in probability theory for the sequence $(f(X_n))$ in the
\begin{theorem}\label{theoreme_proba}
Under the same assumptions than in theorem~\ref{theoreme_equidistribution}.

Denote, for any continuous function $f$ on the torus, $\overline f= f-\int f\di\nu$ and, for any sequence $\underline{x}=(X_n) \in (\T^d)^\N$,
\[
S_n f(\underline{x}) = \sum_{k=0}^{n-1}\overline{f}(X_k) = \sum_{k=0}^{n-1} f(X_k) - n \int f\di\nu
\]
Moreover, for any $t\in [0,1]$, set
\[
\xi_n(t) = \frac 1 {\sqrt {n}} \left(S_if(\underline{x})+ n\left(t-\frac i n\right)f(X_{i})\right) \text{ for }\frac i n \leqslant t\leqslant\frac{i+1} n \text{ and }0 \leqslant i \leqslant n-1
\]
Then, for any continuous function $f$ on the torus and any $x\in \T^d$,
\[
\frac {S_n f(\underline{x})} n \xrightarrow\, 0 \;\;\prob_x-\text{a.e.}
\]
Moreover, for any $\alpha\in]0,1]$ there is $t\in \R_+^\ast$ such that for any $\varepsilon\in ]0,1]$ there is a constant $C$ such that for any $\alpha-$hölder-continuous function $f$ on the torus, any $x\in \T^d$ and any $n\in \N$,
\[
\prob_x \left( \left\{ \underline{x} \in \X^\N\middle|  \left|S_n f(\underline{x})\right|  >n\varepsilon \|f\|_\alpha \right\}\right) \leqslant Ce^{-t \varepsilon^2 n}
\]
Finally, set
\[
\sigma^2(f) := \int_{\T^d} g^2 - (Pg)^2 \di\nu
\]
and then,
\begin{enumerate}
\item\label{item:TCLPSF}If $\sigma^2(f) \not=0$ then for any bounded continuous $F: \mathcal C^0([0,1]) \to \R$ and any $x\in\T^d$,
\[
\esp_x F(\xi_n) \xrightarrow\, \esp F(W_{\sigma^2}) \text{ and }\frac 1 {\ln n } \sum_{k=1}^n \frac 1 k F(\xi_k) \xrightarrow\, \esp F(W_{\sigma^2}) \; \prob_x- \text{a.e.}
\]
Where $W_{\sigma^2}$ denotes Wiener's measure of variance $\sigma^2$.

And for any continuous function $\varphi$ on $\R$ such that $t^2 \varphi(t)$ is bounded and for any $x\in \T^d$,
\[
\frac 1 {\ln n } \sum_{k=1}^n \frac 1 k \varphi\left(\frac{ S_k f(\underline{x})}{\sqrt k}\right)\xrightarrow\, \esp \varphi(W_{\sigma^2}(1)) \;  \prob_x- \text{a.e.}
\]
\item\label{item:variance_zero} If $\sigma^2(f)= 0$ then for any $x\in \T^d$ and any $n\in \N$, $S_nf \in \mathrm{L}^\infty(\prob_x)$ and
\[
 \|S_nf\|_{\mathrm{L}^\infty(\prob_x)} \leqslant 2C\|f\|_\alpha
\]
\end{enumerate}
\end{theorem}

\begin{remark}
The two convergences of $(F(\xi_n))$ in point~$(\ref{item:TCLPSF})$ are respectively called \emph{functional central limit theorem} (FCLT) and \emph{almost-sure functional central limit theorem} (ASFCLT). There is no obvious link between the convergence in law of $(F(\xi_n))$ and the a.e. convergence of it's logarithmic average (see~\cite{BC01} for a criterion). However, note that we have to take a logarithmic mean because of the arc sine law.
\end{remark}

\begin{remark}
The FCLT and the ASFCLT have many corollaries such as the cental limit theorem and the almost sure central limit theorem (taking $F_\varphi(\xi) = \varphi(\xi(1))$ for any continuous and bounded function $\varphi$ on $\R$), the law of the iterated logarithm (see theorem 2.4 in~\cite{Cha96}), a control of $\frac {\max_{k \in [0,n]}S_kf(\underline{x})}{\sqrt n}$ (taking $F(\xi) :=\sup_{t\in [0,1]} \xi(t)$), or an estimation of $\sigma^2(f)$ (taking $\varphi(x) = x^2$).
\end{remark}

Before we continue, we give an example where there is a non-constant function $f$ such that $\sigma^2(f)=0$.

\begin{example}
Let
\[
A=\left(\begin{array}{cc} 2 & 1 \\ 1 & 1\end{array}\right)\text{ and }B =\left(\begin{array}{cc} 0 & 1 \\ -1 & 0\end{array}\right)
\]
Then, the subgroup spanned by $A$ and $B$ is strongly irreducible and proximal.

Let $b_0\in \T^1\setminus\Q/\Z$ be a diophantine number\footnote{There are $C,L\in \R_+^\ast$ such that for any $q\in \N^\ast$, $d(qb_0,0)\geqslant Cq^{-L}$.}, $b=(b_0,0)$ and $\mu = \frac 1 2 \delta_{(A,b)} + \frac 1 2 \delta_{(BA,Bb)}$.

Then, according to proposition~\ref{proposition:critere_diophantien}, the measure $\mu$ satisfies to the assumptions of theorem~\ref{theoreme_equidistribution}.

Let $g$ be the function defined for any $x\in \T^2$ by $g(x) = d(x,0)$. We made everything so that for any $x\in \T^2$, $g(Bx) = g(x)$.

Then, for any $x\in \T^2$,
\[
Pg(x) = \frac 1 2 g(Ax+b) + \frac 1 2 g(BA x+Bb) = g(Ax+b)
\]
And,
\[
\int_\X |Pg(x)|^2 \di\nu(x) = \int_\X |g(Ax+b)|^2 \di\nu(x) = \int_\X |g(x)|^2 \di\nu(x)
\]
Moreover, if we set $f=g-Pg$, then, we just saw that $\sigma^2(f) = \int g^2 -(Pg)^2\di\nu=0$ and for any $x\in \X$, $n\in \N$ and any $(g_1, \dots g_n) \in \{(A,b), (BA,Bb)\}^n$, we have that
\[
g(g_{n+1} \dots g_1 x) =g(A g_n \dots g_1x+b)
\]
And so,
\begin{flalign*}
\sum_{k=0}^{n-1} f(g_k \dots g_1 x) &= g(x) - g(g_n \dots g_1 x) + \sum_{k=0}^{n-1} g(g_{k+1} \dots g_1 x) - g(A g_k \dots g_1x+b) \\&= g(x) - g(g_n \dots g_1 x)
\end{flalign*}
This proves that for any $x\in \X$, the sequence $(\sum_{k=0}^{n-1} f(g_k \dots g_1 x))$ is bounded in $\mathrm{L}^\infty (\prob_x)$.
\end{example}

The results in section~3 of~\cite{Boytore} actually prove that this example is really general.

\bigskip
We will see in sub-section~\ref{soussection:demonstration_theoreme_proba} that theorem~\ref{theoreme_proba} is a quite general corollary of theorem~\ref{theoreme_equidistribution} since we can easily study functions $f$ on the torus that writes $f=g-Pg+\int f\di\nu$ with $g$ continuous and theorem~\ref{theoreme_equidistribution} precisely says that any holder-continuous function can be written in this way.

\medskip
Therefore, the main point of this article is the proof of theorem~\ref{theoreme_equidistribution}. To do so, we use the same method as Bourgain, Furmann, Lindenstrauss and Mozes. In section~\ref{section_non_equidistribution}, we prove that the only obstacle in the equidistribution of the measure $\mu^{\ast n} \ast \vartheta$ is the lower regularity of $\vartheta$ i.e. the existence of points $x$ such that for some $r$ depending on $n$,
\[
\vartheta(B(x,r))\geqslant r^\varepsilon
\]
In particular, if $\mu^{\ast n+m} \ast \vartheta$ is far from Lebesgue's measure then there has to be points $x$ such that
\[
\mu^{\ast m}\ast \vartheta(B(x,r)) \geqslant r^\varepsilon
\]
Then, our assumptions that $\mu$ satisfies an effective shadowing lemma and that $\supp\mu$ is not a subset of $\mathrm{SL}_d(\Z) \ltimes \Q^d/\Z^d$ will allow us to prove in section~\ref{section_stabilisateurs} that this cannot happen when $ r\ll e^{-m} \ll r^\varepsilon$.

The precise proof of the theorem is in subsection~\ref{soussection:demonstration_theoreme_equidistribution}.

\medskip
Finally, in section~\ref{section:fermeture}, we prove proposition~\ref{proposition:critere_diophantien} that is a criterion that shows that under some diophantine conditions on the translations in it's support, a measure satisfies an effective shadowing lemma and we will use this criterion to produce examples of such measures.  
 
\medskip
In the appendix, we state results on the products of random matrices in the case where the action is not irreducible and that we use in section~\ref{section_stabilisateurs}.

\subsection{Some kind of diophantine assumption is necessary}

We already said (and we will prove in section~\ref{section:fermeture}) that a way to guarantee that a measure satisfies an affective shadowing lemma is to require diophantine conditions on the coefficients of the translations of it's support. In this sub-section, we prove that this kind of assumptions is indeed necessary to get theorem~\ref{theoreme_equidistribution}.

\begin{proposition}
Let $a,b\in \mathrm{SL}_d(\Z)$ and $v\in \T^d$. Set $\mu= \frac 1 2 \delta_{(a,0)} + \frac 12 \delta_{(b,v)}$.

Assume that for some $\alpha\in ]0,1]$, there are $C,t\in \R_+^\ast$ such that for any $\alpha-$hölder-continuous function $f$ on the torus and any $n\in \N$,
\[
\sup_{x\in \T^d} \left|P^n f(x) - \int f\di\nu\right| \leqslant Ce^{-tn} \|f\|_\alpha
\]
Then, there are constants $C_0,L\in \R_+^\ast$ such that for any rational point $\frac pq \in \Q^d/\Z^d$,
\[
d\left(v,\frac p q \right)\geqslant \frac{C_0}{q^L}
\]
\end{proposition}

\begin{proof}
For $q\in \N^\ast$ and $x\in \T^d$, we set $X_q = \frac 1 q \Z^d/\Z^d$ and
\[
f_q(x) = 1- \min\left(1, q^{2\alpha}d(x,X_q)^\alpha \right) 
\]
This function is chosen so that it takes the value $1$ on $\frac 1 q \Z^d/\Z^d$, it vanishes on the complementary of the $\frac 1 {q^2}-$neighborhood of $\frac 1 q \Z^d/\Z^d$ and it is hölder-continuous with $\|f_q\|_\alpha \leqslant q^{2\alpha}$.

In particular, we have that, for some constant $C$ depending only on $d$ (and on the distance on $\T^d$),
\[
\int f_q\di\nu \leqslant \sum_{\frac p q \in \frac 1 q \Z^d/\Z^d} \nu\left(B\left(\frac pq, \frac 1 {q^2}\right) \right) \leqslant \frac C{q^{d}}
\]
Moreover, for $\mu^{\otimes \N}-$a.e. $((a_n,b_n))$, we have that
\begin{align*}
\left|f_q\left(\sum_{k=1}^n a_n \dots a_{k+1} b_k\right) - 1 \right|&\leqslant \|f_q\|_\alpha d\left(\sum_{k=1}^n a_n \dots a_{k+1} b_k, X_q \right)^\alpha \\&\leqslant \frac {e^{\alpha Mn}}{(e^M-1)^\alpha} d(v,X_q)^\alpha q^{2\alpha}
\end{align*}
where we noted $e^M = \max(\|A\|,\|B\|)$. Indeed, for any $\frac pq\in X_q$, we have that $f_q(p/q) = 1$ and $\sum_{k=1}^n a_n \dots a_{k+1} b_k$ can be written $Dv$ where $D$ is a matrix with integer coefficients and $\|D\| \leqslant \sum_{k=1}^n \|a_n \dots a_{k+1}\| \leqslant \sum_{k=1}^n e^{M(n-k)}$.
This proves that for any $n\in \N$,
\[
\left| P^n f_q(0) - 1 \right| = \left|\int_{\G} f_q(b) \di\mu^{\ast n}(a,b) - 1 \right|\leqslant \frac {e^{\alpha Mn}}{(e^M-1)^\alpha} d(v,X_q)^\alpha q^{2\alpha}
\]
But, by assumption, we also have that
\[
\left|P^n f_q(0) - \int f_q \di\nu\right| \leqslant Ce^{-tn} \|f_q\|_\alpha \leqslant Ce^{-tn}q^{2\alpha}
\]
So, this proves that for any $n,q\in \N^\ast$,
\[
1 - \frac C{q^d} -Ce^{-tn} q^{2\alpha} \leqslant \frac {e^{\alpha Mn}}{(e^M-1)^\alpha} d(v,X_q)^\alpha q^{2\alpha}
\]
Thus, for any $p\in \Z^d$, any $q\in \N^\ast$ such that $q^d \geqslant 4C$ and any $n$ such that $Ce^{-tn} \leqslant q^{-2\alpha}/4$, we have that
\[
d\left(v,\frac p q\right) \geqslant \frac {e^M - 1}{2^{1/\alpha}e^{Mn} q^{2}}
\]
In particular, for $n = \lfloor \frac 1 t \ln(4Cq^{2\alpha}) \rfloor+1$, we find that, for some constant $C'$ depending only on $M,\alpha,t,C$,
\[
d\left( v, \frac pq\right) \geqslant \frac {C'} {q^{2+2\alpha M/t }}
\]
And this is what we intended to prove.
%
 \end{proof}
 
 \begin{remark}
We can prove the same kind of results for rates more general than $Ce^{-tn}$ and this shows that even convergences slower than exponential require some kind of diophantine assumption.
 \end{remark}

\subsection{Proof of theorem~\ref{theoreme_equidistribution} given the results of sections~\ref{section_non_equidistribution} and~\ref{section_stabilisateurs}} \label{soussection:demonstration_theoreme_equidistribution}

Let $\alpha\in ]0,1]$ and $\varepsilon\in \R_+^\ast$.

According to theorem~\ref{theoreme:BFLMreformule} there are constants $c_0,\varepsilon' \in \R_+^\ast$ with $\varepsilon'<\varepsilon$ such that for any $\vartheta \in \cal M^1(\T^d)$, any $t\in ]0,1]$ and any $n\in \N$ with $n\geqslant c_0(1+\left|\ln t\right|)$,
\[
\mathcal{W}_\alpha\left(\mu^{\ast n} \ast \vartheta, \nu\right) \geqslant t \Rightarrow \vartheta\left( \left\{ x \in \T^d \middle| \vartheta(B(x,r))\geqslant r^{\varepsilon} \right\}\right) \geqslant t^{c_0}
\]
where
\[
r=e^{-(\Lambda_1 +\varepsilon')n} \left( \frac{t}{16} \right)^{1/\alpha}
\]
In particular, for any $m,n\in \N$, any $\vartheta\in\mathcal{M}^1(\T^d)$ any $t\in \R_+^\ast$ small enough and any $C$ large enough,
\[
\mathcal{W}_\alpha\left(\mu^{\ast n+m} \ast \vartheta, \nu\right) \geqslant Ce^{-t n} \Rightarrow \mu^{\ast m}\ast\vartheta\left( \left\{ x \in \T^d \middle| \mu^{\ast m}\ast\vartheta(B(x,r))\geqslant r^{\varepsilon} \right\}\right) \geqslant \left(Ce^{-t n}\right)^{c_0}
\]
for
\[
r=\frac{e^{-(\Lambda_1 +\varepsilon'- \delta /\alpha)n} }{16 ^{1/\alpha}}
\]
But, since the measure satsfies to an effective shadowing lemma, according to proposition~\ref{proposition:mesure_stabilisateur}, there are $C_1,C_2,t_0,L \in \R_+^\ast$ such that for any $x,y\in \T^d$, any $m\in \N$ and any $r\in \R_+^\ast$ with $r\leqslant C_1e^{-Lm}$,
\[
\mu^{\ast m} \left(\left\{ g\in \G\middle| d(gx,y) \leqslant r \right\}\right) \leqslant C_2 e^{-t_0 m}
\]
And so, to get a contradiction, we only need to assume that
\[
r=\frac{e^{-(\Lambda_1 +\varepsilon'- \delta /\alpha)n} }{16 ^{1/\alpha}} \leqslant C_1 e^{-Lm} \text{ and }r^\varepsilon = \frac{e^{-\varepsilon(\Lambda_1 +\varepsilon'- \delta /\alpha)n} }{16 ^{\varepsilon/\alpha}} \geqslant C_2 e^{-t_0m}
\]
And this is always possible for $n=Km$ with $K\in \N$ large enough and $\varepsilon$ small enough.
%
%
%
%

We just proved that there are $C,t\in \R_+^\ast$ and $K\in \N^\ast$ such that for any borelian probability measure $\vartheta$ on $\T^d$ and any $m\in \N$,
\[
\mathcal{W}_\alpha\left( \mu^{\ast (K+1)m} \ast \vartheta,\nu\right) \leqslant Ce^{-tm}
\]
Let $n\in \N$ and let $m,L\in \N$ be such that $n=(K+1)m+L$ and $0\leqslant L< K+1$. Then,
\begin{align*}
\mathcal{W}_\alpha\left( \mu^{\ast n} \ast \vartheta,\nu\right) &\leqslant \mathcal{W}_\alpha\left( \mu^{\ast (K+1)m} \ast \mu^{\ast L}\ast\vartheta,\nu\right) \leqslant C e^{-tm} = Ce^{-\frac t {K+1} (n-L)} \\
&\leqslant C e^t e^{-tn/(K+1)}
\end{align*}
And this finishes the proof of the first part of the theorem.

\medskip
In particular, with $\vartheta= \delta_x$ for some $x\in \T^d$, we get that for any $\alpha-$hölder-continuous function $f$ on $\T^d$ and any $x\in \T^d$,
\[
\left|P^n f(x) - \int f\di\nu\right| \leqslant Ce^{-tn} \|f\|_\alpha
\]

Let $f$ be an $\alpha-$hölder-continuous function on $\T^d$. Set, for any $n\in \N$,
\[
g_n=\sum_{k=0}^{n-1} P^k f-\int f\di\nu
\]
Then,
\[
(I_d-P)g_n = f-\int f\di\nu - \left(P^n f -\int f\di\nu \right)
\]
And so,
\[
\lim_n g_n - Pg_n = f-\int f\di\nu
\]
Moreover, the series is normally convergent since
\[
\sum_n  \left\| P^nf - \int f\di\nu\right\|_\infty \leqslant \frac{C}{1-e^{-t}} \|f\|_\alpha
\]
And so, the function $g=\lim_n g_n$ exists, is continuous and satisfies
\[
g-Pg=f-\int f\di\nu\text{ and }\|g\|_\infty \leqslant  \frac{C}{1-e^{-t}} \|f\|_\alpha
\]

Now, let $\vartheta$ be a $P-$invariant borelian probability measure on $\T^d$. Then, for any hölder-continuous function $f$,
\[
 \int f\di \vartheta = \int P^n f \di\vartheta \xrightarrow[n\to +\infty]\, \int f\di\nu
\]
Where we first used the $P-$invariance of $\vartheta$ and then the dominated convergence theorem since for any $x\in \T^d$, $\lim_n P^n f(x) = \int f\di\nu$ according to the first part of the proof. And, finally, as the hölder-continuous functions are dense in the space of continuous functions on the torus, this proves that $\vartheta = \nu$ and so, $\nu$ is the unique $P-$invariant borelian probability measure on $\T^d$. 

\subsection{Proof of theorem~\ref{theoreme_proba}} \label{soussection:demonstration_theoreme_proba}

\begin{proof}[Proof of the law of large numbers]
This result is a consequence of the uniqueness of the $P-$invariant borelian probability measure seen in theorem~\ref{theoreme_equidistribution}. Indeed, if we manage to prove that for any $x$ and $\prob_x-$a.e. $\underline{x}=(X_n) \in (\T^d)^\N$, the accumulation points of $\nu_{n,{\underline{x}}}:=\frac 1 n \sum_{k=0}^{n-1} \delta_{X_k}$ are $P-$invariant, we will get that they have to be the Lebesgue's measure and so, for any continuous function $f$ on the torus,
\[
\frac 1 n \sum_{k=0}^{n-1} f(X_k) \xrightarrow\, \int f\di\nu \;\;\prob_x-\text{a.e.}
\]
For any continuous function $f$ on the torus, we can compute,
\begin{align*}
\int f\di\nu_{n,{\underline{x}}} - \int Pf\di\nu_{n,{\underline{x}}} &= \frac 1 n \sum_{k=0}^{n-1} f(X_k) - \frac 1 n \sum_{k=0}^{n-1} Pf(X_k) \\
&= \frac 1 n \sum_{k=0}^{n-1} f(X_{k+1}) - Pf(X_k)  + \frac 1n (f(X_0) - f(X_n))
\end{align*}
But, $M_n = \sum_{k=0}^{n-1} f(X_{k+1}) - Pf(X_k)$ is a martingale with bounded increments so $\frac 1 n M_n \xrightarrow\, 0$ a.e. and as $f$ is bounded, we also have that $\frac 1 n (f(X_0) - f(X_n)) \xrightarrow\, 0$ in $\mathrm{L}^\infty(\prob_x)$.

Thus, we just proved that for any $x\in \T^d$ and any continuous function $f$ on $\T^d$, there is $X_f \subset (\T^d)^\N$ such that $\prob_x(X_f) = 1$ and for any $\underline{x}\in X_f$,
\[
\lim_n \int f\di\nu_{n,{\underline{x}}} - \int Pf\di\nu_{n,{\underline{x}}}=0
\]
Let $(f_i)$ be a dense sequence in $\mathcal{C}^0(\T^d)$ and $X_\infty = \cap_i X_{f_i}$. Then, $\prob_x(X_\infty)=1$ and for any $\underline{x}\in X_\infty$ and any $i\in \N$,
\[
\lim_n \int f_i\di\nu_{n,{\underline{x}}} - \int Pf_i\di\nu_{n,{\underline{x}}}=0
\]
So, as the sequence $(f_i)$ is dense, we get that for any continuous function $f$ on $\T^d$ and any $\underline{x} \in X_\infty$,
\[
\lim_n \int f\di\nu_{n,{\underline{x}}} - \int Pf\di\nu_{n,{\underline{x}}} = 0
\]
This proves that for any $\underline{x}\in X_\infty$, the accumulation points of $(\nu_{n,{\underline{x}}})$ are $P-$invariant and so they are equal to $\nu$ and this proves the law of large numbers.
\end{proof}

To prove the remaing part of the theorem, we are going to use Gordin's method and deduce the non-concentration inequality, the FCLT and the ASFCLT from these results for martingales. Indeed, according to theorem~\ref{theoreme_equidistribution}, for any $\alpha\in ]0,1]$, there is a constant $C$ such that for any $\alpha-$hölder-continuous function $f$ on the torus there is a continuous function $g$ such that
\[
f-\int f\di\nu = g-Pg \text{ and }\|g\|_\infty \leqslant C\|f\|_\alpha
\]
Set, for $\underline{x}= (X_n) \in \X^\N$,
\[
S_n f(\underline{x}) = \sum_{k=0}^{n-1} f(X_k) -n\int f\di\nu \text{ and }M_n = \sum_{k=0}^{n-1} g(X_{k+1}) - Pg(X_k)
\]
Then,
\[
S_n f(\underline{x}) = M_n + g(X_0) - g(X_n)
\]
And $M_n$ is a martingale with bounded increments.

\begin{proof}[Proof of the non-concentration inequality]
For any $n\in \N$, we have that
\[
|M_n| \geqslant |S_n f(\underline{x})| - 2 \|g\|_\infty \geqslant |S_n f(\underline{x})| - 2C\|f\|_\alpha
\]
So, using Azuma-Hoeffding's inequality, if $n\varepsilon>2C$, we get that
\begin{flalign*}
I_n(x) :&=\prob_x \left( \left|S_nf(\underline{x}) \right| > n\varepsilon \|f\|_\alpha \right)\leqslant \prob_x\left( \left| M_n \right| \geqslant (n\varepsilon -2C) \|f\|_\alpha\right) \\
& \leqslant 2\exp\left(\frac{ - (n \varepsilon - 2 C)^2\|f\|_\alpha^2}{2n ( 2C\|f\|_\alpha)^2} \right) = 2 \exp\left(-\frac{ n \varepsilon^2}{8 C^2} + \frac{ \varepsilon }{4 C } - \frac 1 {2n}\right) 
\end{flalign*}
And this finishes the proof of this point.\end{proof}

\begin{proof}[Proof of points~\ref{item:TCLPSF}  and~\ref{item:variance_zero}]
As the function $g$ is bounded, the sequence $(S_nf(\underline{x}) - M_n)$ is bounded in $\mathrm{L}^\infty(\prob_x)$ and so it is clear that to prove the FCLT and the AEFCLT, it is enough to study the martingale $M_n$ (that has bounded increments). But, according to the functional central limit theorem for martingales (see corollary 4.1 in~\cite{HH80}) and it's almost sure extension (see~\cite{Cha96}), it is enough to prove the a.e. convergence of the variance (when the limit doesn-t vanish). But, for any $n\in \N^\ast$,
\begin{align*}
\frac 1 n \sum_{k=0}^{n-1} \esp_x\left[\left| M_{k+1} - M_k\right|^2 \middle|X_0, \dots ,X_k \right] &= \frac 1 n \sum_{k=0}^{n-1} \esp_x\left[\left| g(X_{k+1}) - Pg(X_k)\right|^2 \middle|X_0, \dots ,X_k \right] \\& = \frac 1 n \sum_{k=0}^{n-1} P(g^2)(X_k) - (Pg(X_k))^2 
\end{align*}
So, according to the law of large numbers that we already proved and applied to the continuous function $P(g^2) - (Pg)^2$,
\[
\frac 1 n \sum_{k=0}^{n-1} \esp_x\left[\left| M_{k+1} - M_k\right|^2 \middle|X_0, \dots ,X_k \right]  \xrightarrow \,\sigma^2(f):=\int g^2 - (Pg)^2 \di\nu \;\prob_x-\text{a.e.} 
\]
(We used the $P-$invariance of $\nu$ to get that $\int P(g^2) \di\nu = \int g^2 \di\nu$).

And this proves point~\ref{item:TCLPSF} since we suppose in it that $\sigma^2(f) \not=0$.

\medskip
To conclude, remark that, using the $\G-$invariance of $\nu$, we can compute
\begin{align*}
\int_{\G}\int_{\X}  \left| g(\gamma x) - P g(x)\right|^2\di\nu(x)  \di\mu(\gamma) &=  \int_\G  \int_\X g(\gamma x)^2 + g(x)^2 -2Pg(x) g(\gamma x) \di\nu(x) \di\mu(\gamma) \\
&= 2 \int_\X g^2 - (Pg)^2 \di\nu = 2\sigma^2(f)
\end{align*}
And so, if $\sigma^2(f) = 0$, then, as $g$ is continuous, we get that for any $\gamma \in \supp \mu$ and any $x\in \T^d$, $g(\gamma x) = Pg(x)$. This proves that for any $n\in \N$, $M_n=0$ $\prob_x-$a.e. and so, $S_n f(\underline{x}) = g(X_0) - g(X_n)$.
Thus, for any $x\in \T^d$, $S_nf \in \mathrm{L}^{\infty}(\prob_x)$ and
\[
\sup_{x\in \T^d} \sup_{n\in \N} \|S_nf\|_{\mathrm{L}^\infty(\prob_x)} \leqslant 2 \|g\|_\infty \leqslant 2C \|f\|_\alpha
\]
This inequality finishes the proof of point~\ref{item:variance_zero}.
\end{proof}

\section{The non-equidistribution comes from the lower regularity of the measure}\label{section_non_equidistribution}

Like Bourgain, Furmann, Lindenstrauss and Mozes did for the linear random walk on the torus, we are going to prove in this section that if the measure $\mu^{\ast n} \ast \vartheta$ is far from being equidistributed, it is only because of atoms i.e. of points $x\in \T^d$ such that
\[
\vartheta(B(x,r))\geqslant r^\varepsilon
\]
for some $r\in \R_+^\ast$ depending on $n$.

\medskip
More specifically, the aim of this section is to prove the
\begin{theorem}\label{theoreme:BFLMreformule}
Let $\mu$ be a borelian probability measure on $\mathrm{SL}_d(\Z) \ltimes \T^d$. Denote by $\mu_0$ the projection of $\mu$ on $\mathrm{SL}_d(\Z)$ and assume that $\mu_0$ is strongly irreducible, proximal and has an exponential moment. Let $\lambda_1 \in \R_+^\ast$ be the largest Lyapunov exponent of $\mu_0$ (see appendix~\ref{annexe_produits_matrices}).

Then for any $\alpha\in \left]0,1\right]$ and any $\varepsilon \in \R_+^\ast$, there is $c_0,\varepsilon'\in \R_+^\ast$ with $\varepsilon'< \varepsilon$ such that for any $\vartheta \in \cal M^1(\T^d)$, any $t\in \left]0,1\right]$ and any $n\in \N$ with $n\geqslant c_0(1+\left|\ln t\right|)$,
\[
\mathcal{W}_\alpha\left(\mu^{\ast n} \ast \vartheta, \nu\right) \geqslant t \Rightarrow  \vartheta\left( \left\{ x \in \T^d \middle| \vartheta(B(x,r))\geqslant r^{\varepsilon} \right\}\right) \geqslant t^{c_0}
\]
Where we set
\[
r=e^{-(\lambda_1 +\varepsilon')n} \left( \frac{t}{16} \right)^{1/\alpha}
\]
\end{theorem}

In the case of the linear random walk, this statement is a reformulation of an intermediate result (propositions~7.1 and 7.2) of~\cite{BFLM11}. We could prove it for the affine random walk just like they do for the linear one that is to say, by studying, for any borelian probability measure $\vartheta$ on $\T^d$, the set of Fourier-coefficients of $\mu^{\ast n}\ast \vartheta$ and  by remarking that for any $c\in \Z^d$,
\begin{align*}
\widehat{\mu^{\ast n} \ast \vartheta}(c) &= \int_{\T^d}\int_{\mathrm{SL}_d(\Z) \ltimes \T^d} e^{2i\pi \langle c,ax+b \rangle} \di\mu^{\ast n}(a,b) \di\nu(x) \\
&= \int_{\mathrm{SL}_d(\Z) \ltimes \T^d} e^{2i\pi \langle c,b\rangle} \widehat \vartheta(^t ac) \di\mu^{\ast n}(a,b)
\end{align*}
And so,
\[
\left|\widehat{\mu^{\ast n} \ast \vartheta}(c) \right| \leqslant \int_{\mathrm{SL}_d(\Z) \ltimes \T^d} \left|\widehat \vartheta(^t ac)\right| \di\mu^{\ast n}(a,b) = \int_{\mathrm{SL}_d(\Z)}  \left|\widehat \vartheta(^t ac)\right| \di\mu_0^{\ast n}(a)
\]
Where we recall that we denoted by $\mu_0$ the projection of $\mu$ onto $\mathrm{SL}_d(\Z)$.

Thus, if $\left|\widehat{\mu^{\ast n} \ast \vartheta}(c) \right| \geqslant t$, then for many $a$, we also have that $\left|\widehat \vartheta(^t ac)\right|\geqslant t$ and this is the key remark in the proof of BFLM.

\medskip
Instead, we are going to see that this result can also be obtained as a corollary of the one of BFLM for the linear walk : at first, we are going to prove that their result gives informations on the spectral radius of the operator $P$ in $\mathrm{L}^p(\T^d,\nu)$ (even for the affine random walk) and then, that this implies the theorem.

%
%
\subsection{Spectral gap in \texorpdfstring{$\mathrm{L}^p(\T^d)$}{LpTd}}

Let $\G$ be a second countable locally compact group acting measurably on a standard borelian space $\X$ endowed with a $\G-$invariant probability measure $\nu$.

Let $\mu$ be a borelian probability measure on $\G$ and $P$ the Markov operator associated to $\mu$. This is the operator defined for any non-negative borelian function $f$ on $\X$ and any $x\in \X$ by
\[
Pf(x) = \int_\G f(gx) \di\mu(g)
\]
As $\nu$ is a $\G-$invariant probability measure, it is clear that for any $p\in[1,+\infty]$, $1\in \mathrm{L}^p(\X,\nu)$ and $P1=1$. Moreover, we can prove, using Jensen's inequality that $\|P\|_p= 1$. So, we note, for any $p\in ]1,+\infty]$,
\[
L^p_0(\X,\nu):= \left\{ f\in \mathrm{L}^p(\X,\nu)\middle| \int f\di\nu=0 \right\}
\]
and $\rho_p$ the spectral radius of $P$ in $\mathrm{L}^p_0(\X,\nu)$. We say that $P$ has a spectral gap in $\mathrm{L}^p_0(\X,\nu)$ (or, by abuse of notations in $\mathrm{L}^p(\X,\nu)$) if $\rho_p<1$.

In the sequel, we will need a more flexible tool than the spectral gap. This is why, for any $P-$invariant subspace $\Hb$ of $\mathrm{L}^p(\X,\nu)$ endowed with a norm $\|\,.\,\|_\Hb$ such that $P$ is continuous on $ (\Hb,\|\,.\,\|_\Hb)$ and the injection of $(\Hb,\|\,.\,\|_\Hb)$ into $(\mathrm{L}^p, \|\,.\,\|_p)$ is also continuous, we set
\[
\kappa(\mu,\Hb, \mathrm{L}^p\left(\X,\nu\right) ) := -\ln \limsup_{n\to +\infty} \sup_{f\in \Hb\setminus\{0\}} \left( \frac{\|P^n f\|_{p}}{\|f\|_\Hb} \right)^{1/n}
\]
\begin{remark}
The sequence $\left( \sup_{f\in \Hb \setminus\{0\}} \frac{\|P^n f\|_{p}}{\|f\|_\Hb} \right)$ is not sub-multiplicative in general so it may converge to $0$ only at polynomial rate and in this case, we would have that $\kappa(\mu,\Hb, \mathrm{L}^p\left(\X,\nu\right) )=0$. This is impossible if $\Hb = \mathrm{L}^p(\X)$ because in this case, if it converges to $0$, it has to be at exponential rate.
\end{remark}

With this definition, if $(\Hb,\|\,.\,\|_\Hb)= (\mathrm{L}^p(\X,\nu), \|\,.\,\|_p)$, then $e^{-\kappa(\mu,\Hb, \mathrm{L}^p\left(\X,\nu\right) )} = \rho_p$ and, for any $(\Hb,\|\,.\,\|_\Hb)$, we have, since the inclusion of $\Hb$ into $\mathrm{L}^p(\X,\nu)$ is supposed to be continuous,
\[
\kappa(\mu,\Hb, \mathrm{L}^p\left(\X,\nu\right) ) \geqslant -\ln \rho_p
\]

\bigskip
In particular, when $\Hb$ is a subset of $\mathrm{L}^\infty(\X,\nu)$, we can define and study the function $\left(p\mapsto \kappa(\mu,\Hb, \mathrm{L}^p\left(\X,\nu\right) )\right)$.

\begin{remark}\label{remark:kappa_decreasing}
Remind that, according to Hölder's inequality, for any $1\leqslant p\leqslant p'$ and any function $f\in \mathrm{L}^\infty (\X,\nu)$ with $\|f\|_\infty\leqslant 1$,
\[
\|f\|_{p'}^{p'} = \int_{\X} |f|^{p'} \di\nu \leqslant \int_\X |f|^p \di\nu =\|f\|_p^p\text{ and }\|f\|_p \leqslant \|f\|_{p'}
\]
So, we get that the function $\left( p\mapsto \kappa(\mu,\Hb, \mathrm{L}^p\left(\X,\nu\right) )\right)$ is decreasing whereas the function $\left( p\mapsto p\kappa(\mu,\Hb, \mathrm{L}^p\left(\X,\nu\right) )\right)$ is non-decreasing.
\end{remark}

\bigskip
In the same way that we defined $\mathrm{L}^p_0(\T^d)$, we set
\[
\mathcal{C}^{0,\alpha}_0(\T^d) := \left\{f\in \mathcal{C}^{0,\alpha}(\T^d) \middle| \int f\di\nu=0\right\}
\]

The definition of the function $\kappa$ is made to get the
\begin{proposition} \label{proposition:spectral_gap_bflm}
Let $\mu$ be a strongly irreducible and proximal probability measure on $\mathrm{SL}_d(\Z)$ having an exponential moment.

Then, for any $\alpha \in ]0,1]$ small enough,
\[
\lim_{p\to+\infty}  p\, \kappa\left(\mu,\mathcal{C}^{0,\alpha}_0\left(\T^d\right), \mathrm{L}^p\left(\T^d,\nu\right) \right) = \lambda_1 d
\]
%
%
\end{proposition}

This theorem implies in particular that for any $\varepsilon \in \R_+^\ast$, there are $p\in \N$ and $C\in \R_+$ such that for any $n\in \N$ and any $f\in \mathcal{C}^{0,\alpha}_0(\T^d)$,
\[
\|P^n f\|_{\mathrm{L}^p(\T^d)}\leqslant C e^{-(\lambda_1 d - \varepsilon) n/p} \|f\|_\alpha
\]

\begin{proof}
First of all, since $\mu$ has an exponential moment, for any $\alpha\in ]0,1]$ small enough, any $f\in \mathcal{C}^{0,\alpha}(\T^d)$ and any $x\in \T^d$,
\[
|Pf(x)| \leqslant \|f\|_\infty
\]
and for any $y\in \T^d$,
\begin{align*}
\left|Pf(x) - Pf(y)\right| &\leqslant \int_\G |f(gx) - f(gy)| \di\mu(g) \leqslant \|f\|_\alpha \int_\G d(gx,gy)^\alpha \di\mu(g)\\
& \leqslant \|f\|_\alpha d(x,y)^\alpha \int_\G \|g\|^\alpha \di\mu(g)
\end{align*}
So, $Pf$ is $\alpha-$hölder-continuous and $P$ is a continuous operator on $\mathcal{C}^{0,\alpha}_0(\T^d)$.

\medskip
According to the result of Bourgain, Furmann, Lindenstrauss and Mozes in~\cite{BFLM11} (that we use as stated in proposition 4.5 in~\cite{Boytore}), we have that for any $\alpha\in ]0,1]$ and any $\varepsilon \in \R_+^\ast$, there is a constant $C$ such that for any $n\in \N$, any $t\in ]0,1]$ with $n\geqslant -C\ln t$ and any $f\in \mathcal{C}^{0,\alpha}(\T^d)$ with $\int f\di\nu=0$,
\[
\left\{x\middle||P^n f(x)|\geqslant t \|f\|_\alpha \right\} \subset \bigcup_{\substack{\frac p q \in \Q^d/\Z^d\\q\leqslant Ct^{-C}}} B\left(\frac pq, e^{-(\lambda_1 - \varepsilon)n} \right) 
\]
In particular, for any $L\in \N^\ast$,
\begin{align*}
\int |P^n f|^L \di\nu& \leqslant (t\|f\|_\alpha)^L + \nu\left( \left\{x\middle||P^n f(x)|\geqslant t \|f\|_\alpha \right\}\right) \|f\|_\infty^L\\
& \leqslant \left(t^L +(Ct^{-C})^d e^{-(\lambda_1 - \varepsilon)dn} \right)\|f\|_\alpha^L 
\end{align*}
And so, taking $t=e^{-\delta n}$ with $\delta\in \R_+^\ast$ small enough and $L\in \N$ large enough, we find that for some constant $C$,
\[
\int|P^n f|^L \di\nu \leqslant Ce^{-(\lambda_1 - 2\varepsilon)dn} \|f\|_\alpha^L
\]
And this proves (reminding that the limit exists according to remark~\ref{remark:kappa_decreasing}) that
\[
\lim_p  p\, \kappa\left(\mu,\mathcal{C}^{0,\alpha}_0\left(\T^d\right), \mathrm{L}^p\left(\T^d,\nu\right) \right) \geqslant \lambda_1 d
\]
We are now going to prove the other inequality. Let $\delta \in ]0,1/4]$, $f\in \mathcal{C}^\infty(\T^d)$ such that $f=1$ on $B(0,\delta)$, $\|f\|_\infty \leqslant 1$ and $\int f =0$.

Then, for any $\varepsilon \in \R_+^\ast$ and any $x\in B(0, e^{-(\lambda_1 +\varepsilon)n} \delta)$, we have that
\begin{align*}
P^n f(x) &= \int_\G \un_{\|g\|\leqslant e^{(\lambda_1 + \varepsilon)n}} f(gx) \di\mu^{\ast n}(g)+ \int_\G \un_{\|g\|\geqslant e^{(\lambda_1 + \varepsilon)n}} f(gx) \di\mu^{\ast n}(g) \\
&= \mu^{\ast n}\left(\left\{ g\middle| \|g\|\leqslant e^{(\lambda_1 + \varepsilon)n} \right\}\right) +  \int_\G \un_{\|g\|\geqslant e^{(\lambda_1 + \varepsilon)n}} f(gx) \di\mu^{\ast n}(g) \\
& \geqslant 1 - 2 \mu^{\ast n}\left(\left\{ g\middle| \|g\|\geqslant e^{(\lambda_1 + \varepsilon)n} \right\}\right)
\end{align*}
But, according to theorem~\ref{theoreme:LDP_produits_reductible}, there are $C,t\in \R_+^\ast$ such that
\[
\mu^{\ast n}\left(\left\{ g\middle| \|g\|\geqslant e^{(\lambda_1 + \varepsilon)n} \right\}\right) \leqslant Ce^{-tn}
\]
And so, for $n\in\N$ large enough, we have that for any $x\in \T^d$,
\[
|P^n f(x)| \geqslant \left(1-2Ce^{-tn}\right) \un_{B(0,e^{-(\lambda_1 + \varepsilon)n}\delta)}(x)
\]
In particular, for any $L\in \N$,
\[
\int |P^n f(x)|^L \di\nu \geqslant \left( 1-2Ce^{-tn} \right)^L \nu(B(0, e^{-(\lambda_1 + \varepsilon)n}\delta)) =  \left( 1-2Ce^{-tn} \right)^L e^{-(\lambda_1 + \varepsilon)dn} \delta^d
\]
And this proves that
\[
\lim_p  p\, \kappa\left(\mu,\mathcal{C}^{0,\alpha}_0\left(\T^d\right), \mathrm{L}^p\left(\T^d,\nu\right) \right) \leqslant \lambda_1 d
\]
And this finishes the proof of the proposition.
\end{proof}

We are now going to extend the previous result to measures on $\mathrm{SL}_d(\Z) \ltimes \T^d$ by proving the
\begin{corollary}\label{coro:ts}
Let $\mu$ be a borelian probability measure on $\G:=\mathrm{SL}_d(\Z) \ltimes \T^d$. Let $\mu_0$ be the projection of $\mu$ onto $\G_0:=\mathrm{SL}_d(\Z)$ and assume that $\mu_0$ is strongly irreducible, proximal and has an exponential moment.
%

Then, for any $\alpha \in ]0,1]$ small enough,
\[
\lim_{p\to+\infty}  p\, \kappa\left(\mu,\mathcal{C}^{0,\alpha}_0\left(\T^d\right), \mathrm{L}^p\left(\T^d,\nu\right) \right) \geqslant \lambda_1 d
\]
\end{corollary}

\begin{remark}
Theorem~\ref{theoreme_equidistribution} actually proves that for any measure $\mu$ satisfying it's assumptions we have that for any $f\in \mathcal{C}^{0,\alpha}(\T^d)$ with$\int f\di\nu=0$, any $n\in \N$ and any $p\in \N^\ast$,
\[
\int |P^n f|^p \di\nu \leqslant C^p e^{-tpn} \|f\|_\alpha^p
\]
And so, for any $p\in [1,+\infty[$,
\[
\kappa\left( \mu,\mathcal{C}^{0,\alpha}_0\left(\T^d\right), \mathrm{L}^p\left(\T^d,\nu\right) \right) \geqslant t
\]
And, in particular,
\[
\lim_{p\to +\infty} p\;\kappa\left(\mu,\mathcal{C}^{0,\alpha}_0\left(\T^d\right), \mathrm{L}^p\left(\T^d,\nu\right) \right)  = +\infty
\]
\end{remark}

We are going to prove this result in three steps. First, we are going to prove it for trigonometric functions, then, for regular ones and last, for hölder-continuous functions.

\begin{lemma}
Let $\mu$ be a borelian probability measure on $\G$ and $\mu_0$ it's projection on $\G_0$. Denote by $P$ the Markov operator associated to $\mu$ and by  $P_0$ the one associated to $\mu_0$.

Then, for any $c \in \Z^d$, any $n\in \N$ and any $L\in \N$,
\[
\int_{\T^d} |P^n e_c|^{2L} \di\nu \leqslant \int_{\T^d} |P_0^n e_c|^{2L} \di\nu
\]
Where, for $c\in \Z^d$, $e_c$ is the function defined for $x\in \T^d$ by
\[
e_c(x) := e^{2i\pi \langle c,x\rangle}
\]
\end{lemma}

\begin{proof}
Using Fubini's theorem, we can make the following computation
\begin{align*}
\int_\X &|P^n_0 e_c(x)|^{2L} \di\nu (x)=\int_\X (P^n_0 e_c(x))^L \overline{(P^n_0f(x))^L}\di\nu(x) \\
&= \int_\X \int_{\G^{2L}}e_c(a_1x) \dots e_c(a_L x) \overline{ e_c(a_{L+1}x) \dots e_c(a_{2L}x)}  \di\mu_0^{\ast n}(a_1) \dots \di\mu_0^{\ast n}(a_{2L})\di\nu(x) \\
&=  \int_{\G^{2L}} \int_\X e^{2i\pi \langle c,\left(a_1+ \dots + a_L - (a_{L+1} + \dots + a_{2L})\right)x\rangle} \di\nu(x) \di\mu_0^{\ast n}(a_1) \dots \di\mu_0^{\ast n}(a_{2L}) \\
&= \int_{\G^{2L}} \un_{\{^t(a_1+ \dots + a_L - (a_{L+1} + \dots + a_{2L}))c=0\}} \di\mu_0^{\ast n}(a_1) \dots \di\mu_0^{\ast n}(a_{2L}) \\
\end{align*}
Doing the same kind of computations for the measure $\mu$, and noting, to simplify notations, for $(a_1,b_1), \dots, (a_{2L}, b_{2L}) \in \supp\mu$, $a_i^j = a_i +\dots + a_j$ and $b_i^j = b_i + \dots + b_j$, we find that
\begin{align*}
\int_\X |P^n e_c|^{2L} \di\nu &= \int_{\G^{2L}} e^{2i\pi\langle c,b_1^L - b_{L+1}^{2L}\rangle}\un_{\{^t(a_1^L - a_{L+1}^{2L})c=0\}} \di\mu^{\ast n}(a_1,b_1) \dots \di\mu^{\ast n}(a_{2L},b_{2L}) \\
& \leqslant \int_{\G^{2L}} \un_{\{^t(a_1^L- a_{L+1}^{2L})c=0\}} \di\mu_0^{\ast n}(g_1) \dots \di\mu_0^{\ast n}(g_{2L}) =\int_\X |P^n_0 e_c|^{2L} \di\nu
\end{align*}
Where the last inequality comes from the first part of the proof and precisely gives what we intended to prove.
\end{proof}

For $s\in \R_+^\ast$, we denote by $\mathcal{H}^s(\T^d)$ the Sobolev space of exponent $s$.

\begin{lemma}\label{lemme:tsregulier}
With the same assumptions than in corollary~\ref{coro:ts}, for any $s\in \R_+^\ast$ large enough and any $\varepsilon \in \R_+^\ast$, there is $L\in \R_+$ such that for any $f\in \mathcal{H}^s(\T^d)$ and any $n\in \N$,
\[
\int \left| P^n f-\int f\di\nu\right|^{2L} \di\nu \leqslant C e^{-(\lambda_1d - \varepsilon)n} \|f\|_{\mathcal{H}^s(\T^d)}^{2L}
\]
\end{lemma}

\begin{proof}
Let $f\in \mathcal{H}^s(\T^d)$. Then, by definition, we can expand $f$ in Fourier series : $f = \sum_{c\in \Z^d} \widehat f(c) e_c$ with $\|f\|_{\mathcal{H}^s}:=\left(\sum_{c\in \Z^d} (1+\|c\|^2)^{s/2} |\widehat{f}(c)|^2\right)^{1/2} <+\infty$ and so, for any $L\in \N^\ast$,
\[
\left(\int_{\T^d} |P^n f|^{2L} \di\nu\right)^{1/2L} \leqslant \sum_{c\in \Z^d} |\widehat f(c)| \left(\int_{\T^d } |P^n e_c|^{2L} \di\nu \right)^{1/2L}
\]
Using the previous lemma and proposition~\ref{proposition:spectral_gap_bflm}, we get that for any $\varepsilon \in \R_+^\ast$, there is $L\in \R_+$ such that for any $c\in \Z^d\setminus\{0\}$,
\[
\int_{\T^d} |P^n e_c(x)|^{2L} \di\nu(x) \leqslant C e^{-(\lambda_1d- \varepsilon)n} \|c\|^{2L}
\]
Combining this inequality with the previous one, we get that for any $f \in \cal H^s(\T^d)$ with $\widehat f(0) = \int f\di\nu=0$,
\[
\left(\int_{\T^d} |P^n f|^{2L}\di\nu\right)^{1/2L} \leqslant C^{1/2L} e^{-(\lambda_1 d-\varepsilon)n/2L} \sum_{c\in \Z^d\setminus\{0\}} |\widehat f(c)|\|c\| \leqslant C' e^{-(\lambda_1 d-\varepsilon)n/2L} \|f\|_{\mathcal{H}^s}
\]
For some constant $C'$ depending on $d$, $s$, $\mu$, $L$ but non on $f$.
\end{proof}

\begin{proof}[End of the proof of corollary~\ref{coro:ts}]
According to Jackson-Bernstein's lemma, for any $\alpha\in ]0,1]$ and any $s\in \R_+$ large enough, there is a constant $C$ such that for any $f\in \mathcal{ C}^{0,\alpha}(\T^d)$, there is a sequence $(f_m)\in \mathcal{H}^s(\T^d)^{\N}$ such that for any $m\in \N^\ast$,
\[
\int f\di\nu = \int f_m\di\nu, \quad \|f-f_m\|_\infty \leqslant \frac C{m^\alpha} \|f\|_\alpha \text{ and }\|f_m\|_{\mathcal{H}^s} \leqslant Cm\|f\|_\alpha
\] 
This implies that for any $x\in \T^d$ and any $m,n\in \N^\ast$,
\[
|P^n f_m(x)| \geqslant |P^n f (x)| - \frac C {m^\alpha} \|f\|_\alpha
\]
%
%

Let $\varepsilon \in \R_+^\ast$, $m\in \N^\ast$ and $t=\frac {2C}{m^\alpha}$. Then, using the equality $t-C/m^\alpha = t/2$ and lemma~\ref{lemme:tsregulier}, we get that
\begin{align*}
\int |P^n f|^{2L} \di\nu & \leqslant (t\|f\|_\alpha)^{2L} + \nu\left(\{ |P^n f|\geqslant t\|f\|_\alpha\}\right) \|f\|_\infty^{2L}\\
&\leqslant \left(t^{2L} + \nu\left(\left\{ |P^n f_m| \geqslant \left(t-\frac C {m^\alpha}\right) \|f\|_\alpha \right\}\right) \right) \|f\|_\alpha^{2L}\\
& \leqslant \left(t^{2L} + \left( \frac 2 {t\|f\|_\alpha} \right)^{2M} \int |P^n f_m|^{2M} \di\nu \right)\|f\|_\alpha^{2L} \\
& \leqslant  \left(t^{2L} + \left( \frac 2 {t\|f\|_\alpha} \right)^{2M} C^{2M} m^{2M} e^{-(\lambda_1 d -\varepsilon)n} \|f\|_\alpha^{2M}\right)\|f\|_\alpha^{2L}
\end{align*}
So, for $m=e^{\delta n}$, we get that for some constant $C'$,
\[
\int |P^n f|^{2L} \di\nu \leqslant C' \left(e^{-\delta \alpha 2L n} + e^{ \delta (1+\alpha)2M n -(\lambda_1 d - \varepsilon) n} \right)\|f\|_\alpha^{2L}
\]
And so, for $\delta$ small enough and $L$ large enough, we get that
\[
\int |P^n f|^{2L} \di\nu \leqslant C e^{-(\lambda_1 d-2\varepsilon)n} \|f\|_\alpha^{2L}
\]
And this is what we intended to prove.
\end{proof}

\subsection{Equidistribution, lower regularity and spectral gap.}

In this subsection, we finish the proof of theorem~\ref{theoreme:BFLMreformule} by studying the link between equidistribution and the lower regularity of the measure when the spectral gap is large.

\begin{lemma} \label{lemme:normeholder}
Under the same assumptions as in theorem~\ref{theoreme:BFLMreformule}, for any $\varepsilon \in \R_+^\ast$, there are $C,t \in \R_+^\ast$ such that for any $x,y\in \T^d$, any $f\in \mathcal{C}^{0,\alpha}(\T^d)$ and any $n\in \N$,
\[
|P^n f(x) - P^n f(y)| \leqslant \left( e^{\alpha(\lambda_1 + \varepsilon) n} d(x,y)^\alpha + Ce^{-tn} \right) \|f\|_\alpha
\]
\end{lemma}

\begin{proof}
Let's compute, for any $x,y\in \T^d$, $n\in \N$, $f\in \mathcal{C}^{0,\alpha}(\T^d)$ and $\varepsilon \in \R_+^\ast$,
\begin{align*}
\left|P^n f(x) - P^n f(y) \right| &= \left|\int_\G f(gx) - f(gy) \di\mu^{\ast n}(g) \right| \\
&\leqslant \int_{\G} \un_{\|g\|\leqslant e^{(\lambda_1 + \varepsilon)n}} |f(gx) - f(gy)| \di\mu^{\ast n}(g) \\& \retrait\retrait+ \int_\G \un_{\|g\|\geqslant e^{(\lambda_1 + \varepsilon)n}} |f(gx) - f(gy)| \di\mu^{\ast n}(g) \\
& \leqslant m_\alpha(f) \int_\G \un_{\|g\|\leqslant e^{(\lambda_1 + \varepsilon)n}} d(gx,gy)^\alpha \di\mu^{\ast n}(g)\\& \retrait\retrait + 2\|f\|_\infty \mu^{\ast n} \left( \left\{ g\in \G \middle| \|g\|\geqslant e^{(\lambda_1 + \varepsilon)n} \right\}\right) \\
& \leqslant \left( d(x,y)^{\alpha} e^{\alpha(\lambda_1 + \varepsilon)n} + 2  \mu^{\ast n} \left( \left\{ g\in \G \middle| \|g\|\geqslant e^{(\lambda_1 + \varepsilon)n} \right\}\right) \right) \|f\|_\alpha
\end{align*}
Where we used the fact that for any $x,y\in \T^d$ and any $g\in \G$,
\[
d(gx,gy) \leqslant \|g\| d(x,y)
\]
To conclude, we use theorem~\ref{theoreme:LDP_produits_reductible} and we get that there are $C,t\in \R_+^\ast$ such that
\[
\mu^{\ast n} \left( \left\{ g\in \G \middle| \|g\|\geqslant e^{(\lambda_1 + \varepsilon)n} \right\}\right) \leqslant Ce^{-tn} \qedhere
\]
\end{proof}

We are now ready to prove theorem~\ref{theoreme:BFLMreformule}.

The idea of the proof is that if we have some point $x_0$ of the torus such that $|P^nf(x_0)| \geqslant t$, then, on a neighborhood $B(x_0,r)$  for $r\approx e^{-\lambda_1 n}$ we also have that $|P^n f(x)| \approx t$. But, the control on $\kappa(\mu,\mathcal{C}^{0,\alpha}, \mathrm{L}^p)$ implies that $\nu(x| |P^n f(x)| \geqslant t) \approx e^{-\lambda_1 dn}$ and so, we have that $\nu(x| |P^n f(x)| \geqslant t) \approx e^{-\lambda_1 dn}\approx \nu(B(x_0,r))$ and this proves that $\{ x ||P^n f(x)| \geqslant t\}$ cannot be much bigger than $B(x_0,r)$.

\begin{proof}[Proof of theorem~\ref{theoreme:BFLMreformule}]
Let $\alpha,t\in ]0,1]$, $\vartheta \in \mathcal{M}^1(\T^d)$ and $n\in \N$. As, for any $0<\alpha'<\alpha$ the inclusion of $\mathcal{C}^{0,\alpha}(\T^d)$ into $\mathcal{C}^{0,\alpha'}(\T^d)$ is continuous, we may assume without any loss of generality that $\alpha$ is small enough so that corollary~\ref{coro:ts} holds.

Assume that $\mathcal{W}_\alpha(\mu^{\ast n} \ast \vartheta,\nu)\geqslant t$.

By definition, there is $f\in \mathcal{C}^{0,\alpha}(\T^d)$ with $\|f\|_\alpha\leqslant 1$ and such that
\[
\left|\int P^n f\di\vartheta - \int f \di\nu \right| \geqslant \frac t 2
\]
We can assume without any loss of generality that $\int f\di\nu=0$ and $\|f\|_\alpha \leqslant 2$.
And this proves that
\[
\int_{\T^d} \left|P^n f(x)\right| \di\vartheta(x) \geqslant  \left|\int_{\T^d}P^n f(x) \di\vartheta(x)\right| \geqslant \frac t 2
\]
We set, for any $n\in \N$ and $t\in ]0,1]$,
\[
X_{n,t} := \left\{x\in \T^d \middle| \left|P^n f(x)\right| \geqslant t\right\}
\]
Then, using that $\|P^n f\|_\infty \leqslant \|f\|_\infty \leqslant 2$, we find that
\[
\frac t 2 \leqslant \int_{\X} |P^n f(x)| \di\vartheta(x) \leqslant \frac t 4 + 2 \vartheta\left(X_{n,t/4} \right) 
\]
And so,
\[
 \vartheta\left( X_{n,t/4}  \right) \geqslant \frac t 8
\]
Moreover, according to lemma~\ref{lemme:normeholder}, for any $\varepsilon_2\in \R_+^\ast$, there are $C,t_0\in \R_+^\ast$ such that for any $x\in X_{n,t/4}$ and any $y\in \T^d$, we have that
\[
|P^n f(y)| \geqslant \frac t 4-e^{\alpha(\lambda_1 + \varepsilon_2)n} d(x,y)^{\alpha} - Ce^{-t_0 n} \geqslant \frac t 8 -e^{\alpha(\lambda_1 + \varepsilon_2)n} d(x,y)^{\alpha}
\]
Since we can take $c_0$ so large that $Ce^{-t_0 n} \leqslant \frac t 8$ for $n\geqslant c_0(1+|\ln t|)$.

In particular, noting $r= e^{-(\lambda_1 + \varepsilon_2 )n}\left(\frac t{16} \right)^{1/\alpha}$, we have that for any $x\in X_{n,t/4}$, $B(x, r)\subset X_{n,t/16}$.

Moreover, according to the classical covering results, there is a constant $C(d)$, depending only on $d$ and points $x_1, \dots, x_N \in \T^d$ such that
\[
X_{n,t/4} \subset \bigcup_{i=1}^N B(x_i,r) \subset X_{n,t/16}
\]
and the union has multiplicity at most $C(d)$. 

This implies in particular that
\[
\sum_{i=1}^N \un_{B(x_i,r)} \leqslant C(d) \un_{X_{n,t/16}}
\]
And so, taking the integral against the measure $\nu$ and using the equality $\nu(B(x,r)) = r^d$, we get
\[
Nr^d \leqslant C(d) \nu(X_{n,t/16})
\]
To sum-up, we found points $x_1,\dots, x_N$ with $ N\leqslant\frac{C(d) \nu(X_{n,t/16})}{r^d}$, such that
\[
\vartheta\left(\bigcup_{i=1}^N B(x,r)\right) \geqslant \frac t 8
\]
So, noting $\mathcal{I}:=\{i\in [1,N] | \vartheta(B(x_i,r)) \geqslant \frac t{16N} \}$, we get that
\[
\vartheta\left(\bigcup_{i\in \mathcal{I}} B(x_i,r)\right) \geqslant \frac t {16}
\]
And finally, for any $x\in \bigcup_{i\in \mathcal{I}} B(x_i,r)$, there is, by definition of $\mathcal{I}$, some $i \in \mathcal{I}$ such that
\[
B(x,2r) \supset B(x_i, r)
\]
And so,
\[
\vartheta(B(x,2r)) \geqslant \frac t{16N}
\]
In conclusion, we proved that
\[
\vartheta\left(\left\{x\in \T^d \middle| \vartheta(B(x,2r)) \geqslant \frac t{16 N} \right\}\right) \geqslant \frac t {16}
\]
To finish, note that we have that
\[
N \leqslant C(d) \frac{\nu(X_{n,t/16})}{r^d} = C(d) \left(\frac{16} t\right)^{d/\alpha} e^{(\lambda_1 + \varepsilon_2)dn} \nu(X_{n,t/16})
\]
And that, according to Markov's inequality and corollary~\ref{coro:ts}, for any $\varepsilon_1 \in \R_+^\ast$, there are $C,L\in \R_+$ such that
\begin{align*}
\nu(X_{n,t/16}) &\leqslant \left(\frac {16} t\right)^{2L} \int_\X |P^n f(x)|^{2L} \di\nu(x)\leqslant \left(\frac {16} t\right)^{2L} Ce^{-(\lambda_1d-\varepsilon_1)n} \|f\|_\alpha^{2L} \\& \leqslant\left(\frac {32} t\right)^{2L} Ce^{-(\lambda_1d-\varepsilon_1)n} 
\end{align*}
This proves that for some constant $C$ depending on $\varepsilon_1,d,\mu$,
\[
N \leqslant \frac{ C}{t^C} e^{(\varepsilon_1 + \varepsilon_2 d ) n}
\]
And so, taking $\varepsilon_1,\varepsilon_2$ small enough and $c_0$ large enough, we get that
\[
\frac{t}{16 N} \geqslant \frac{ t^{C+1}}C e^{-(\varepsilon_1+ \varepsilon_2 d) n} \geqslant (2r)^\varepsilon
\]
and this is what we intended to prove.
\end{proof}

\section{Measure of points-stabilizers}\label{section_stabilisateurs}

The aim of this section is to prove the following
\begin{proposition}\label{proposition:mesure_stabilisateur}
with the same assumptions as in theorem~\ref{theoreme_equidistribution}, there are $C,t \in \R_+^\ast$ such that for any $x\in \T^d$ and any $n\in \N$,
\[
\mu^{\ast n} \otimes \mu^{\ast n} \left(\left\{g_1,g_2 \in \G \middle| g_1 x = g_2 x\right\}\right) \leqslant Ce^{-tn}
\]
\end{proposition}

This proposition will be a direct corollary of lemmas~\ref{lemma:mesure_stab_tore} and~\ref{lemma:mesure_stab_Rd} since if $\mu$ is not concentrated on $\mathrm{SL}_d(\Z)\ltimes \Q^d$, then the measure $\mu_1$ of lemma~\ref{lemma:mesure_stab_tore} is not concentrated on $\mathrm{SL}_d(\R) \ltimes \{0\}$.

\medskip
To evaluate the measure of the stabilizer of a point, we are going to lift the situation from $\T^d$ to $\R^d$ since we know better the products of random elements of $\mathrm{SL}_d(\R) \ltimes \R^d$ (through the theory of products of random matrices since this group can be identified to a subgroup of $\mathrm{SL}_{d+1}(\R)$) than those of elements of $\mathrm{SL}_d(\Z) \ltimes \T^d$.

\medskip
To understand what we are going to do, remark that if $\mu = \frac 1 2 \delta_{(g_1,v_1)} + \frac 1 2 \delta_{(g_2,v_2)}$ with $g_1,g_2\in \mathrm{SL}_d(\Z)$, $v_1\in \Q^d/\Z^d$ and $v_2\in\T^d$ such that the coefficients of $v_2$ and $1$ are $\Q-$linearly independent, then, for $\mu^{\ast n}-$a.e. $(g,v)\in \mathrm{SL}_d(\Z)\ltimes \T^d$, we can write
\[
v=M_1 v_1 + M_2 v_2
\]
with $M_i \in \mathcal{M}_d(\Z)$.

So, in particular, noting $\overrightarrow{v_1}, \overrightarrow{v_2}$ some representatives of $v_1,v_2$ in $\R^d$, we get that if $v=0$ in $\T^d$, then there is $p\in \Z^d$, such that,
\[
M_1 \overrightarrow{v_1} + M_2 \overrightarrow{v_2} = p
\]
As we assumed that the coefficients of $v_2$ and $1$ are $\Q-$linearly independent and that $v_1\in \Q^d$, we get that $M_2 \overrightarrow{v_2} = 0$ and $M_1 \overrightarrow{v_1} = p$.

So, we set
\[
\mu_1 = \frac 1 2 \delta_{(g_1,0)} + \frac 1 2 \delta_{(g_2,\overrightarrow{ v_2})} \in \mathcal{M}^1 \left(\mathrm{SL}_d(\R) \ltimes \R^d \right)
\]
and what we just proved is that
\[
\mu^{\ast n} \left( \mathrm{SL}_d(\Z) \ltimes \{0\}\right) \leqslant \mu_1^{\ast n} \left( \mathrm{SL}_d(\R) \ltimes \{0\}\right)
\]
And so we are let with a problem on probability measures on $\mathrm{SL}_d(\R) \ltimes \R^d$.

\medskip
Thus, to lift the situation from $\T^d$ to $\R^d$, we are going to project the translation part of elements in the support of $\mu$ onto a complementary subspace of $\Q^d$ in the $\Q-$vector space $\R^d$. To do so, we fix some $\Q-$linear projection $\pi_\Q:\R^d \to \R^d$ onto $\Q^d$ and we remark that for any $v\in \R^d/\Z^d$, $v-\pi_\Q v \in \R^d$ is well defined and the application $\mathrm{SL}_d(\Z) \ltimes \T^d \ni (g,v) \mapsto (g,v-\pi_\Q v) \in \mathrm{SL}_d(\Z) \ltimes \R^d$ is a group-morphism.

Now, we can prove the following
\begin{lemma} \label{lemma:mesure_stab_tore}
Let $\mu$ be a borelian probability measure on $\mathrm{SL}_d(\Z) \ltimes \T^d$.

Then, for any $n\in \N$ and any $y\in \T^d$,
\[
\mu^{\ast n} \otimes\mu^{\ast n} \left( \left\{(g_1,g_2) \middle| g_2^{-1} g_1 \in \mathrm{Stab}(y)\right\} \right) \leqslant \mu_1^{\ast n} \otimes \mu_1^{\ast n} \left( \left\{(g_1,g_2) \middle| g_2^{-1} g_1 \in \mathrm{Stab}(y-\pi_\Q y)\right\} \right)
\]
where $\mu_1$ is the measure on $\mathrm{SL}_d(\Z) \ltimes \R^d$ defined by $\mu_1 (A) = \mu(\varphi^{-1}(A))$ where $\varphi: \mathrm{SL}_d(\Z) \ltimes \T^d \to \mathrm{SL}_d(\Z) \ltimes \R^d$ is the function defined by
\[
\varphi(g,v) = (g,v-\pi_\Q v)
\]
and $\pi_\Q$ is some $\Q-$linear projection onto $\Q^d$.
\end{lemma}

\begin{proof}
As $\varphi$ is a morphism, we only need to prove that for any $g\in \mathrm{SL}_d(\Z) \ltimes \T^d$ and any $y\in \T^d$, if $gy=y$ then $\varphi(g) (y-\pi_\Q y) = (y-\pi_\Q y)$.

Write $g=(a,b) $ with $a\in \mathrm{SL}_d(\Z)$ and $b\in \T^d$. Then, $\varphi(g) = (a,b-\pi_\Q b)$.
So,
\[
\varphi(g) (y-\pi_\Q y) = a (y-\pi_\Q y) + b-\pi_\Q b
\]
But, $gy=y$, so, noting $\overrightarrow{b}, \overrightarrow{y}$ some representatives of $b,y$, we get that there is $p\in \Z^d$ such that
\[
a \overrightarrow{y} + \overrightarrow{b} = \overrightarrow{y}+p
\]
Projecting onto $\Q^d$, we also get that
\[
a \pi_\Q\overrightarrow{y} + \pi_\Q\overrightarrow{b} = \pi_\Q\overrightarrow{y}+p
\]
And this proves that
\[
a (I_d-\pi_\Q)\overrightarrow{y}+ (I_d-\pi_\Q)\overrightarrow{b} = (I_d-\pi_\Q)\overrightarrow{y}
\]
Finally, as $(I_d-\pi_\Q)\overrightarrow{y}$ and $ (I_d-\pi_\Q)\overrightarrow{b}$ don't depend on the choices of the representatives of $y$ and $b$, this finishes the proof of the fact that
\[
a(y-\pi_\Q y) + b-\pi_\Q b = y-\pi_\Q y
\]
and this finishes the proof of the lemma.
%
%
%
%
%
%
%
\end{proof}

\begin{lemma}\label{lemma:mesure_stab_Rd}
Let $\mu$ be a borelian probability measure on $\mathrm{SL}_d(\R) \ltimes \R^d$ that is not concentrated on $\mathrm{SL}_d(\R) \ltimes\{0\}$ and has an exponential moment. Let $\mu_0$ be the projection of $\mu$ onto $\mathrm{SL}_d(\R)$ and assume that $\mu_0$ is strongly irreducible and proximal.

Then, there are $C,t \in \R_+^\ast$ such that for any $x\in \R^d$ and any $n\in \N$,
\[
\mu^{\ast n} \otimes \mu^{\ast n}\left( \left\{(g_1,g_2) \middle| g_1x= g_2 x\right\} \right) \leqslant Ce^{-tn}
\]
\end{lemma}

\begin{proof}
We denote by $\tilde \mu$ the measure on $\mathrm{SL}_d(\R)\ltimes \R^d$ defined by $\tilde\mu(A) = \mu(A^{-1})$ for any borelian subset $A$ of $\mathrm{SL}_d(\R)\ltimes \R^d$ and where $A^{-1} := \{g^{-1 } | g\in A\}$.

Let $\lambda_1 \geqslant \dots \geqslant \lambda_d$ be the Lyapunov exponents of $\mu_0$ (see appendix~\ref{annexe_produits_matrices}).

Then, the largest Lyapunov exponent of $\tilde \mu$ is $-\lambda_d$. Moreover, as $\lambda_1 \geqslant \dots \geqslant \lambda_d$, $\lambda_1>0$ and $\lambda_1 + \dots + \lambda_d = 0$, we have that $\lambda_d <0$ and so, $-\lambda_d>0$.

%

Let $n\in \N$, $\varepsilon\in \R_+^\ast$ and $x\in \R^d$. We can compute,
\begin{align*}
I_n(x):&=\mu^{\ast n} \otimes \mu^{\ast n} \left(\left\{g_1,g_2 \in \G \middle| g_1 x=g_2x \right\}\right) = \mu^{\ast n} \otimes \tilde\mu^{\ast n} \left(\left\{(g_1,g_2) \in \G^2 \middle| g_2g_1 x=x \right\}\right) \\
&\leqslant \mu^{\ast n} \otimes\tilde\mu^{\ast n} \left(\left\{(g_1,g_2)\middle| \frac{1+\|g_2^{-1} g_1x\|}{1+\|x\|} \leqslant e^{(\lambda_1 - \lambda_d- \varepsilon)n}\right\} \right) 
\end{align*}
But
\[
\frac{1+\|g_2 g_1x\|}{1+\|x\|} = \frac{1+\|g_2 g_1x\|}{1+\|g_1x\|}\frac{1+\|g_1x\|}{1+\|x\|}
\]
So, we obtain that
\begin{align*}
I_n(x) &\leqslant \mu^{\ast n}\left(\left\{ g_1\middle| \frac{1+\|g_1x\|}{1+\|x\|} \leqslant e^{(\lambda_1 - \varepsilon/2)n} \right\}\right) \\& \retrait\retrait+ \int_\G \tilde\mu^{\ast n}\left(\left\{ g_2\middle| \frac{1+\|g_2g_1x\|}{1+\|g_1x\|} \leqslant e^{(-\lambda_d - \varepsilon/2)n} \right\}\right) \di\mu^{\ast n}(g_1)
\end{align*}
And we can conclude with corollary~\ref{corollaire:longueurs_translations} applied to the measures $\mu$ and $\tilde\mu$.
\end{proof}

\section{Effective shadowing lemmas} \label{section:fermeture}

The aim of this section is to prove a criterion to produce measures satisfying an effective shadowing lemma.

First of all, we recall the
\begin{definition}
Let $\mu$ be a borelian probability measure on $\G$.

We say that $\mu$ satisfies an effective shadowing lemma if for ant $C',t' \in \R_+^\ast$, there are $C_1,C_2,M,t,L \in \R_+^\ast$ such that for any $x,y\in \T^d$, any $r\in \R_+^\ast$ and any $n\in \N$ with $r\leqslant C_1e^{-Ln}$, if
\[
\mu^{\ast n} \left(\left\{ g\in \G\middle| d(gx,y) \leqslant r \right\}\right) \geqslant C_2 e^{-t n}
\]
then there are $x',y'\in \T^d$ such that $d(x,x'),d(y,y')\leqslant re^{Mn}$ and
\[
\mu^{\ast n} \left(\left\{ g\in \G\middle| gx'=y'\right\}\right)\geqslant C'e^{-t'n}
\]
\end{definition}

The criterion that we are going to prove will use the diophantine properties of the translation parts of the elements in $\supp\mu$. More specifically, we want a condition that ensures that if $(g_1,v_1), (g_2,v_2)\in \supp\mu^{\ast n}$ are such that $v_1$ and $v_2$ are close (in some sense) then $v_1=v_2$. This is why we make the following 

\begin{definition}[Diophantine subsets of $\T^d$]\label{definition:mesures_diophantiennes}
Let $d\in \N^\ast$ and $B\subset \T^d$ a finite subset.

We that that $B$ is $(C,L)-$diophantine if for any non zero $(M_b)\in \Z^B$,
\[
d\left(\sum_{b\in B} M_b b,0\right)\leqslant \frac{ C}{\max_b |M_b|^{L}} \Rightarrow \sum_b M_b b=0
\]
%
More generally, we say that $B$ is diophantine if it is $(C,L)-$diophantine for some $C,L\in \R_+^\ast$.
\end{definition}

\begin{remark}
With this definition, a diophantine subset can contain rational points.
\end{remark}

\begin{example}\label{exemple:element_diophantien}
Let $b\in \T^1 \setminus\Q/\Z$. Asking $\{b\}$ to be diophantine is asking for $C,L\in \R_+^\ast$ such that for any $q\in \N^\ast$,
\[
d(qb,0) >\frac{ C}{q^L}
\]
The name comes from this property.
\end{example}

\begin{remark}
Let $d\geqslant 2$ and $B$ a subset of $\T^d$. It is not the same thing to say that $B$ is diophantine and that $\{\text{coefficients of }b| b\in B\}$ is a diophantine subset of $\T^1$ (consider $B=\{(b_1,b_2)\}$ with $b_1$ diophantine and $b_2$ not). This last property is stronger but this is the one that we will need in the sequel and we refer to lemma~\ref{lemma:equivalence_defi_diophantiennes} for more details.
\end{remark}

\begin{example} \label{example:diophantien_generique}
Let $N\in \N^\ast$. Then, for a.e. $b_1, \dots b_N\in \T^1$, the set $\{b_1, \dots b_N\}$ is diophantine.
\end{example}

We are now ready to state the main result of this section
\begin{proposition} \label{proposition:critere_diophantien}
Let $\mu$ be a borelian probability measure on $\mathrm{SL}_d(\Z) \ltimes \T^d$ and let $\mu_0$ be it's projection on $\mathrm{SL}_d(\Z)$. Assume that $\mu_0$ is strongly irreducible, proximal and that it has an exponential moment and that $\{\text{coefficients of }b| (a,b)\in \supp\mu\}$ is a diophantine subset of $\T^1$.

Then, $\mu$ satisfies an effective shadowing lemma.
\end{proposition}

To prove this proposition, we first come back to the difference for a subset $B$ of $\T^d$ between being diophantine and having elements whose coefficients form a diophantine subset of $\T^1$.

\begin{lemma} \label{lemma:equivalence_defi_diophantiennes}
Let $B $ be a finite subset of $ \T^d$ and
\[
F:= \{ \text{coefficients of }B\} \subset \T^1
\]
Then, $F$ is $(C,L)-$diophantine if and only if there is $C'\in \R_+^\ast$ such that for any non zero $(M_b) \in \mathcal{M}_d(\Z)^B$,
\[
d\left(\sum_{b\in B} M_b b,0\right) \leqslant \frac{C'}{(\max_b \|M_b\|)^L} \Rightarrow \sum_b M_bb=0
\]
\end{lemma}

\begin{proof}
First, assume that $F$ is $(C,L)-$diophantine and set $C'= C/(d|B|)^L$.

Let $0\not=(M_b) \in \mathcal{M}_d(\Z)^B$ be such that
\[
d(\sum_{b\in B} M_b b,0) \leqslant \frac{C'}{(\max_b \|M_b\|)^L}
\]
Each coefficient of $\sum_b M_bb$ is a sum of elements of $F$ multiplied by integers that are smaller than $d |B| \max_b \|M_b\|$. In other words, for any coefficient of $\sum_b M_bb$, we get a sum $\sum_{f\in F} L_ff$ with $|L_f| \leqslant d|B| \max_b \|M_b\|$ and
\[
d\left( \sum_f L_ff,0\right) \leqslant \frac {C'}{\max_b \|M_b\|^L} \leqslant \frac{C'(d|B|)^L}{\max_{f}|L_f|^L} = \frac{C}{(\max_f |L_f|)^L}
\]
And as $F$ is $(C,L)-$diophantine, this implies that $\sum_f L_ff=0$ and so, as this is true for any coefficient of $\sum_b M_bb$, we get that $\sum_b M_bb=0$.

Reciprocally, if there is $C' \in \R_+^\ast$ such that for any $0\not=(M_b) \in \mathcal{M}_d(\Z)^B$ 
\[
d\left(\sum_{b\in B} M_bb,0\right) \leqslant \frac{C'}{(\max_b \|M_b\|)^L} \Rightarrow \sum_b M_bb=0
\]
Then, we set $C= C'/d^L$ and let $0\not=(L_f) \in \Z^F$ be such that
\[
d\left(\sum_{f\in F}L_ff,0\right) \leqslant \frac{ C}{(\max_f |L_f|)^L}
\]
For any element $f$ of $F$, choose an element $(b(f),i(f))$ in 
\[
\{(b,i) \in B\times [1,d] | f\text{ is the }i-\text{th coefficient of }b \}
\]
Now, for any $b\in B$, we denote by $M_b$ the matrix where we set $L_f$ in the $i-$th column if there is $f\in F$ such that $(b(f),i(f)) = (b,i)$ and $0$ otherwise.

Thus, by definition,
\[
\sum_{b\in B} M_bb = \left(\begin{array}{c}
\sum_f L_ff\\ \vdots \\\sum_f L_ff
\end{array} \right)
\]
And $\max_b \|M_b\|\leqslant d \max_f |L_f|$ so
\[
d\left(\sum_b M_bb,0\right) \leqslant \frac{ C}{(\max_f|L_f|)^L} \leqslant \frac{ C'}{(\max_b\|M_b\|)^L}
\]
And this proves, as $B$ is $(C,L)-$diophantine, that $\sum_b M_bb=0$ and so we also get that $\sum_f L_f f=0$ and this finishes the proof of the lemma.
\end{proof}

From now on, we set
\begin{equation}\label{equation:BmuFmu}
B(\mu):=\{b| (a,b)\in \supp\mu\} \text{ and }F(\mu):=\{ \text{coefficients of }b|b\in B(\mu)\}
\end{equation}

To prove proposition~\ref{proposition:critere_diophantien}, we are going to use some control of the translation parts of elements of $\supp\mu^{\ast n}$. To do so, for any $Q\in \N^\ast$ and any finite subset $B$ of $  \T^d$, we set
\[
X_Q(B) = \left\{ \frac {p + \sum_{b\in B} M_bb}q\middle| p\in\Z^d,\;q\in \Z,\;|q|\leqslant Q,\;(M_b)\in \mathcal{M}_d(\Z)^B,\;\max_{b\in B}\|M_b\|\leqslant Q\right\}
\]

Thus some $x\in \T^d$ belongs to $X_Q(B) $ if there is $q\in \N^\ast$ with $q\leqslant Q$ such that each coefficient of $qx$ can be obtained from translations by coefficients of elements of $B$ with multiplicities smaller than $Q$.

This definition is made so that for any $M$ large enough we have that for any $n\in \N^\ast$, with large probability, any element $(g,v)\in \supp\mu^{\ast n}$ is such that $v\in X_{e^{Mn}}(B(\mu))$.

To make this idea more precise, we set, for any $M\in \R_+^\ast$ and $n\in \N$,
\[
\G_n^M := \left\{( a_i) \in \G_0^n \middle|\max_{k\in [1,n]} \max_{1\leqslant i_1<\dots<i_k \leqslant n} \|a_{i_k}\dots a_ {i_1}\| \leqslant e^{Mn}\right\}
\]
The aim of the following lemma is to prove that elements of $\G_n^M$ are generic.
\begin{lemma} \label{lemma:controle_translations}
Let $\mu$ be a borelian probability measure on $\G$ and let $\mu_0$ be the image of $\mu$ on $\mathrm{SL}_d(\Z)$. Assume that $\mu_0$ has an exponential moment.

Then, for any $M\in \R_+^\ast$ large enough, there is $ t\in \R_+^\ast$ such that for any $n\in \N$,
\[
\mu^{\otimes n}_0\left(\G_n^M\right) \geqslant 1 - e^{-tn}
\]
Moreover, if $B(\mu)$ is finite (see equation~\eqref{equation:BmuFmu}), then, for any $(a_1,b_1), \dots ,(a_n,b_n) \in \supp\mu$ such that $(a_i) \in \G_n^M$,
\[
\sum_{k=1}^n a_n \dots a_{k+1} b_k \in X_{e^{(M+1)n}}(B(\mu))
\]
\end{lemma}

\begin{proof}
First, remark that as $\mu_0$ is concentrated on $\mathrm{SL}_d(\Z)$, for $\mu_0-$a.e. $g\in \G_0$, $\|g\|\geqslant 1$ and so, noting $\delta\in \R_+^\ast$ such that $\int_{\G_0} \|g\|^\delta \di\mu_0(g)<+\infty$, we have that for any $M\in \R_+^\ast$ and any $n\in \N$,
\begin{align*}
\mu_0^{\otimes n} (\G_n^M):&=\mu_0^{\otimes n} \left(\left\{ (g_i) \in \G_0^\N\middle| \max_{k\in [1,n]} \max_{1\leqslant i_1<\dots<i_k \leqslant n} \|g_{i_k}\dots g_ {i_1}\| \geqslant e^{Mn}\right\}\right) \\
&  \leqslant \mu_0^{\otimes n}\left( \left\{ (g_i) \in \G_0^\N\middle|  \|g_n\|\dots \|g_1\|\geqslant e^{Mn}\right\}\right) \\
& \leqslant e^{-\delta Mn} \left( \int_{\G_0} \|g\|^\delta \di\mu_0(g)\right)^n
\end{align*}
And this finishes the proof of the first part of the lemma.

To prove the second one, take $(a_1,b_1), \dots ,(a_n,b_n) \in \supp\mu$ such that $(a_i) \in \G_n^M$.

Then, we can write
\[
\sum_{k=1}^{n} a_n \dots a_{k+1} b_k = \sum_{b'\in B(\mu)} \left(\sum_{k=1}^n \un_{b_k=b'} a_n \dots a_k \right)b'
\]
and, for any $b'\in B(\mu)$,
\[
\left\|\sum_{k=1}^n \un_{b_k=b'} a_n \dots a_k \right\| \leqslant ne^{Mn}
\]
This proves that $\sum_{k=1}^n a_n \dots a_{k+1} b_k \in X_{ne^{Mn}}(B(\mu))$.
\end{proof}

In the sequel, we will have to control sums of the form $\sum c_i a_i$ where the $a_i$ are $\mu^{\ast n}-$generic and $c_i \in \Z$. To do so, we will use the
\begin{lemma}[\cite{BFLM11}] \label{lemma:BFLM}
Let $\mu_0$ be a strongly irreducible and proximal borelian probability measure on $\mathrm{SL}_d(\Z)$ having an exponential moment.

Set
\[
\G_0^1 := \left\{ a\in \G_0 \middle|  \left|\frac 1 n \ln \|a\|- \lambda_1\right|\leqslant \varepsilon \right\}
\]
and
\[
\G_0^d := \left\{ (a_i)\in \G_0^d \middle| \forall i\in [1,d],\;a_i\in \G_0^1\text{ and } \det\left(\sum_i c_i a_i\right)\not=0 \right\}
\]
where we put $c_1 = 1-d$ and $c_i=1$ for $i\in [2,d]$.

Then, for any $\varepsilon \in \R_+^\ast$, there are $C,t\in \R_+^\ast$ such that for any $n\in \N$,
\[
\mu^{\ast n}(\G_0^1) \geqslant 1-Ce^{-tn}
\text{ et }\mu^{\ast n} \otimes \dots \otimes \mu^{\ast n}(\G_0^d) \geqslant  1 - Ce^{-tn}
\]
\end{lemma}

\begin{proof}
In~\cite{BFLM11}, this lemma is not stated as we do it here but it corresponds to the lemmas~4.3, 4.6 and 7.9 and to some part of the proof of proposition 7.3.
\end{proof}

We are now ready to prove proposition~\ref{proposition:critere_diophantien}. The proof consists in two lemmas. In lemma~\ref{lemma:approx}, we study points that are far from points of $X_{e^n}(B(\mu))$ and in lemma~\ref{lemma:diophantien} we study points that are close to it.

\begin{lemma}\label{lemma:approx}
With the assumptions of proposition~\ref{proposition:critere_diophantien}, for any $M\in \R_+^\ast$ large enough there are $C,t\in \R_+^\ast$ such that for any $r\in ]0,1]$, any $n\in \N$, any $x,y\in \T^d$, if
\[
\mu^{\ast n}\left(\left\{g\in \G \middle| d(gx,y)\leqslant r\right\}\right) \geqslant Ce^{-tn}
\]
then there is $x'\in X_{e^{Mn}}(B(\mu))$ such that
\[
d(x,x') \leqslant re^{Mn}
\]
\end{lemma}

\begin{proof}
The idea that makes the demonstration work is that if $a_1 x + b_1, \dots , a_d x+b_d$ are close to each-others, then setting $c_1=1-d$ and $c_i=1$ for $i\in [2,d]$, we get that $\sum_i c_i a_i x$ is close to $\sum_i c_i b_i$. But, according to lemma~\ref{lemma:BFLM}, with large probability, the matrix $\sum_i c_i a_i$ is invertible and, according to lemma~\ref{lemma:controle_translations}, the $b_i$ belong to $X_{e^{Mn}}(B(\mu))$, and so $x\in X_{e^{Mn}}(B(\mu))$.

We keep the notations $\G_0^1$ and $\G_0^d$ from lemma~\ref{lemma:BFLM} and the notation $\G_n^M$ from lemma~\ref{lemma:controle_translations}. Then, for any $\varepsilon \in \R_+^\ast$ and any large enough $M$ we denote by $C,t$ the constants given by lemma~\ref{lemma:BFLM} and by $t_0$ the one given by lemma~\ref{lemma:controle_translations}.

Then, we set
\[
\G_\ast := \left\{ ((a_1,b_1), \dots, (a_d,b_d))\in \G^d \middle|(a_i)\in \G_0^d \text{ and for any }i,\; b_i \in X_{e^{(M+1)n}}(B(\mu)) \right\}
\]
By definition, we have that
\[
\mu^{\ast n} \otimes  \dots \otimes \mu^{\ast n}(\G_\ast) \leqslant Ce^{-tn} + d e^{-t_0n}
\]
So, we can compute, for any $x,y\in \T^d$,
\begin{align*}
I_n(x,y):&=\left(\mu^{\ast n}\left(\left\{g\in \G \middle| d(gx,y)\leqslant r\right\}\right)\right)^d \\ &= \int_{\G^d}\un_{B(y,r)}(a_1 x+b_1) \dots \un_{B(y,r)} (a_dx+b_d) \di\mu^{\ast n}(a_1,b_1) \dots \di\mu^{\ast n}(a_d,b_d)  \\
& \leqslant Ce^{-tn} + de^{-t_0n} + \int_{\G_\ast} \prod_{i=1}^d\un_{B(y,r)}(a_i x+b_i) \di\mu^{\ast n}(a_1,b_1) \dots \di\mu^{\ast n}(a_d,b_d)
\end{align*}
Thus, if $M$ is large enough, we can find $C,t \in \R_+^\ast$ such that for any $r\in ]0,1]$, any $n\in \N$, if
\[
\mu^{\ast n}\left(\left\{g\in \G \middle| d(gx,y)\leqslant r\right\}\right) \geqslant Ce^{-tn}
\]
then
\[
\int_{\G_\ast} \prod_{i=1}^d\un_{B(y,r)}(a_i x+b_i) \di\mu^{\ast n}(a_1,b_1) \dots \di\mu^{\ast n}(a_d,b_d) >0
\]
In particular, there is $((a_1,b_1), \dots, (a_d,b_d))\in \G^\ast$ such that for any $i$,
\[
d(a_i x+b_i, y) \leqslant r
\]
Now, we let $\overrightarrow{x}, \overrightarrow{y}, \overrightarrow{b_i}$ be representatives of $x,y,b_i$ in $\R^d$. We have that for any $i\in [1,d]$ there is $p_i \in \Z^d$ such that
\[
\|a_i \overrightarrow{x} + \overrightarrow{b_i} - \overrightarrow{y} - p_i\| \leqslant r
\]
So, noting $c_1 = 1-d$ and $c_i = 1$ for $i\not=1$, we get that
\[
\| (\sum_i c_i a_i)x + \sum_i c_i b_i - \sum_i c_ip_i\| \leqslant 2d r
\]
But, by definition of $\G_0^d$, $\det\left(\sum_i c_i a_i \right)\in \Z^\ast$.

So,
\[
1 \leqslant \left|\det\left( \sum_i c_i a_i\right)\right| \leqslant \left\|\sum_i c_i a_i\right\|^d \leqslant d^{2d} e^{(\lambda_1 + \varepsilon)dn}
\]
Let $U = \sum_i c_i a_i$. Then $U$ is invertible and we can write,
\[
U^{-1} = \frac 1 {\det U} V
\]
with $V \in \mathcal{M}_d(\Z)$ and for some constant $C(d)$ depending only on $d$, $\|V\| \leqslant C(d)\|U\|^{d-1}$.

Thus, we get that
\[
\|x+U^{-1} \sum_i c_i b_i - U^{-1}\sum_i c_ip_i \| \leqslant 2dr\|U\|
\]
To conclude, we only need to remark that we can write
\[
b_i = \sum_{b \in B(\mu)} M_g^i b \text{ with }\max_b \|M_b^i\|\leqslant e^{(M+1)n}
\]
And this proves that
\[
-U^{-1}\sum_i c_i b_i +U^{-1} \sum_i c_i p_i =\frac{ \sum_i c_i V p_i}{\det(U)} + \frac{ \sum_{b\in B( \mu)} \left(\sum_i c_i VM_b^i\right)b } {\det(U)}
\]
And that
\begin{align*}
\max_{b\in B(\mu)} \left\| \sum_i c_i V M_b^i\right\| &\leqslant 2d \|V\| \max_{b}\|M_b^i\| \leqslant 2d C(d) e^{(M+1)n} e^{(\lambda_1 + \varepsilon)(d-1)n}
\end{align*}
So, maybe for some $M'\geqslant M$ we get that
\[
x':=-U^{-1}\sum_i c_i b_i +U^{-1} \sum_i c_i p_i  \in X_{e^{M'n}}(\supp\mu)
\]
and $d(x,x') \leqslant re^{M'n}$.
\end{proof}

\begin{lemma}\label{lemma:diophantien}
With the same assumptions as in proposition~\ref{proposition:critere_diophantien}, there are $C_0,L$ such that for any $M\in \R_+$ large enough, there are $C,t\in \R_+^\ast$ such that for any $n\in \N^\ast$, any $r\in ]0,1]$ with
\[
r \leqslant \frac {C_0}{e^{MLn}}
\]
we have that for any $y\in \T^d$, any $x'\in X_{e^{(M+1)n}}(B(\mu))$ and any $x\in \T^d$ with $d(x,x')\leqslant r$,
\[
\left(\mu^{\ast n}\left(\left\{ g\in \G\middle| d(gx,y) \leqslant r\right\}\right) \right)^2\leqslant Ce^{-tn} + \mu^{\ast n}\otimes \mu^{\ast n}\left(\left\{g_1,g_2 \middle| g_1x'=g_2 x'\right\} \right)
\]
\end{lemma}

\begin{proof}To simplify our notations, we set $B= B(\mu)$ and $F=F(\mu)$. 
let $x,x',y\in \T^d$ as in the proposition. By definition of $X_Q(B)$, there is $(M_v)\in \mathcal{M}_d(\Z) ^{B}$ with $\max \|M_v\| \leqslant Q$, there is $p\in \Z^d$ and $q\in \N^\ast$ with $|q|\leqslant Q$ such that
\[
x' = \frac {p+ \sum_{b\in B} M_b b}q
\]

Let, for any $\varepsilon \in \R_+^\ast$ and $M\in \N$,
\[
\G_\ast:= \left\{ (a,b) \in \G\middle| \|a\|\leqslant e^{(\lambda_1 + \varepsilon)n} \text{ et }b\in X_{e^{(M+1)n}}(B)\right\}
\]
Lemmas~\ref{lemma:BFLM} and~\ref{lemma:controle_translations} prove that for any $\varepsilon \in \R_+^\ast$ and any $M\in \R_+$ large enough, there are $C,t\in \R_+^\ast$ such that for any $n\in \N$,
\[
\mu^{\ast n}(\G_\ast) \geqslant 1-Ce^{-tn}
\]
Let's compute
\begin{align*}
I_n(x,y):&=\left(\mu^{\ast n}\left(\left\{ g\in \G\middle| d(gx,y) \leqslant r\right\}\right) \right)^2 \\&= \int_{\G^2} \un_{B(y,r)}(g_1 x) \un_{B(y,r)} (g_2 x) \di\mu^{\ast n}(g_1) \di\mu^{\ast n}(g_2) \\
& \leqslant 2\mu^{\ast n}(\G_\ast) + \int_{\G^2} \un_{\G_\ast}(g_1)\un_{\G_\ast}(g_2) \un_{d(g_1x,g_2x)\leqslant 2r} \di\mu^{\ast n}(g_1) \di\mu^{\ast n}(g_2) \\
& \leqslant 2Ce^{-tn} + \int_{\G^2} \un_{\G_\ast}(g_1)\un_{\G_\ast}(g_2) \un_{d(g_1x,g_2x)\leqslant 2r} \di\mu^{\ast n}(g_1) \di\mu^{\ast n}(g_2)
\end{align*}
Moreover, for $g_1 = (a_1,b_1)$ and $g_2=(a_2,b_2)$, we have that
\[
d(g_1 x',g_2 x')\leqslant d(g_1 x',g_1 x) + d(g_1x,g_2 x) + d(g_2 x,g_2 x') \leqslant \|a_1\| r + d(g_1 x,g_2 x) + \|a_2\|r
\]
And this proves that
\[
I_n(x,y) \leqslant 2 Ce^{-tn} + \int_{\G^2}\un_{\G_\ast}(g_1)\un_{\G_\ast}(g_2) \un_{d(g_1x',g_2x')\leqslant 3re^{(\lambda_1 + \varepsilon)n}} \di\mu^{\ast n}(g_1) \di\mu^{\ast n}(g_2)
\]
To conclude, we only need to prove that, under the diophantine condition, if $x'\in X_{e^{(M+1)n}}$, $g_1,g_2\in \G_\ast$ are such that $d(g_1 x',g_2 x') \leqslant 3re^{(\lambda_1+ \varepsilon)n}$ then $g_1 x'=g_2 x'$.

To do so, remark that if $x'\in X_{e^{(M+1)n}}(B)$ and $g=(a,b)\in \G_\ast$, alors, $a x'+b \in X_{2 e^{(M+1 +\lambda_1 + \varepsilon)n}}(B)$

Thus, we have that $(g_1 - g_2) x' \in X_{4e^{(M+1+\lambda_1 +\varepsilon)n}}(B)$ and that $d((g_1 - g_2)x',0) \leqslant 3re^{(\lambda_1 + \varepsilon)n}$.

In particular, there is $q\in \N^\ast$ with $q\leqslant 4e^{(M+1+\lambda_1 +\varepsilon)n}$ such that $q(g_1-g_2)x'$ is a sum of elements of $b$ multiplied on the left by matrices of norm smaller than $(4e^{(M+1+\lambda_1 +\varepsilon)n})^2$ and
\[
d(q(g_1-g_2)x',0) \leqslant |q|3r e^{(\lambda_1 + \varepsilon)n} \leqslant 12re^{(M+1+2\lambda_1 +2\varepsilon)n}
\]
So, as $F$ is $(C,L)-$diophantine, according to lemma~\ref{lemma:equivalence_defi_diophantiennes}, there are constants $C',L$ such that if
\[
12re^{(M+1+2\lambda_1 +2\varepsilon)n} \leqslant \frac{C'} {(4e^{(M+1+\lambda_1 +\varepsilon)n})^{2L}}
\]
then we have that $qg_1 x'=qg_2 x'$. But, this proves that $g_1 x' = g_2 x' + \frac pq$ for some $p\in \Z^d$. And
\[
d\left(\frac pq,0\right)=d(g_1x',g_2 x') \leqslant 3r e^{(\lambda_1 + \varepsilon)n}
\]
So, if $\frac 1 {|q| }> 3r e^{(\lambda_1 +\varepsilon)n}$, we have that $\frac pq=0$ and so, $g_1 x'=g_2 x'$.
\end{proof}

\begin{proof}[End of the proof of proposition~\ref{proposition:critere_diophantien}]
According to lemma~\ref{lemma:approx}, we have that for any $M\in \R_+^\ast$ large enough, there are $C,t\in \R_+^\ast$ such that for any $r\in ]0,1]$, any $n\in \N$ and any $x,y\in \T^d$, if
\[
\mu^{\ast n}\left(\left\{g\in \G \middle| d(gx,y)\leqslant r\right\}\right) \geqslant Ce^{-tn}
\]
then there is $x'\in X_{e^{Mn}}(B(\mu))$ such that
\[
d(x,x') \leqslant re^{Mn}
\]
But, in this case, we have, according to lemma~\ref{lemma:diophantien}, that if
\[
re^{Mn} \leqslant \frac{C_0}{e^{LMn}},
\]
then
\[
\left(\mu^{\ast n}\left(\left\{ g\in \G\middle| d(gx,y) \leqslant r\right\}\right) \right)^2\leqslant Ce^{-tn} + \mu^{\ast n}\otimes \mu^{\ast n}\left(\left\{g_1,g_2 \middle| g_1x'=g_2 x'\right\} \right)
\]
And this is what we intended to prove since
\[
\mu^{\ast n}\otimes \mu^{\ast n}\left(\left\{g_1,g_2 \middle| g_1x'=g_2 x'\right\} \right) = \int_\G \mu^{\ast n} \left( \left\{ g_1 \middle| g_1 x' = g_2 x' \right\}\right) \di\mu^{\ast n}(g_2)\qedhere
\]
\end{proof}

\appendix
\section{Products of random matrices} \label{annexe_produits_matrices}

In this section, we are going to recall some of the properties of products of random matrices.

To do so, we fix some finite dimensional $\R-$vector space $\V$ that we endow with an euclidian norm.

Let $\mu$ be a borelian probability measure on $\G:=\mathrm{GL}(\V)$.

We set, for any $g\in \G$, $N(g) = \max(\|g\|,\|g^{-1}\|)$ and we say that $\mu$ has a moment of order $1$ if
\[
\int_\G \ln N(g) \di\mu(g)<+\infty
\]
and that it has an exponential moment if there is some $\varepsilon \in \R_+^\ast$ such that
\[
\int_\G N(g)^\varepsilon \di\mu(g) <+\infty
\]
Remark that there is some constant $C$ depending only on $\dim(\V)$ such that for any $g\in \mathrm{SL}_d(\V)$, $\|g^{-1} \| \leqslant C \|g\|^{\dim\V}$ and so, if $\mu$ is a measure on $\mathrm{SL}(\V)$, it is enough to ask that for some $\varepsilon\in \R_+^\ast$,
\[
\int_\G \|g\|^\varepsilon \di\mu(g)<+\infty
\]

We would like to study the product $g_n \dots g_1$ where $(g_i)$ is an iid sequence of law $\mu$.

The first result in this direction is Oseledec's theorem :
\begin{theorem}
Let $\mu$ be a borelian probability measure on $\G:=\mathrm{GL}(\V)$ having a moment of order $1$.

Then, there are $m_1 , \dots , m_r \in \N^\ast $ with $\sum_i m_i=\dim\V$, there is $\Lambda_1>\dots>\Lambda_r \in \R$ and some measurable function from $\G^\N$ into the space of flags $\V:=\V_1 \supset \dots \supset  \V_{r+1}:=\{0\}$ of $\V$ such that $\dim \V_i =  \sum_{j=i}^{r} m_j$ and such that for $\mu^{\otimes \N}-$a.e. $\omega =(g_n) \in \G^\N$,
\begin{itemize}
\item $g_1 \V_i^\omega = \V_i^{\vartheta \omega}$ where $\vartheta$ is the shift on $\G^\N$.
\item For any $x\in\V_i^\omega \setminus \V_{i+1}^\omega$,
\[
\lim_n \frac 1 n \ln\|g_n \dots g_1 x\|=\Lambda_i
\]
\end{itemize}
\end{theorem}

We call \emph{Lyapunov exponents} the paramaters $\Lambda_1, \dots ,\Lambda_r$ and we note
\[
\lambda_1 = \dots = \lambda_{m_1} = \Lambda_1, \quad \lambda_{m_1+1} = \dots = \lambda_{m_1+m_2} = \Lambda_2, \quad\text{ etc.}
\]
However, if $x\in \V$, this theorem doesn't say anything on the behavior of $\|g_n \dots g_1 x\|$ for some generic $(g_n)\in \G^\N$ because we have no information on the sequences $(g_n)$ such that $x\in \V_i^{(g_n)}$. To avoid this problem, we usually assume that the subgroup of $\G$ spanned by the support of $\mu$ acts irreducibly on $\V$ (it doesn't fix any non-trivial subspace of $\V$). And in this case, Furstenberg proved the
\begin{theorem}[\cite{BL85}]
Let $\mu$ be a borelian probability measure on $\G$ having a momen tof order $1$ and whose support generates a group acting irreducibly on $\V$.

Then, for any $x\in \V\setminus\{0\}$,
\[
\frac 1 n \ln \|g_n \dots g_1 x\| \xrightarrow\, \lambda_1 \text{ a.e.}
\]
\end{theorem}

This irreducibility assumption is not good enough for us since we will identify $\mathrm{SL}_d(\R) \ltimes \R^d$ with the subgroup
\[
\left(\begin{array}{cc}
\mathrm{SL}_d(\R) & \R^d \\ 0 & 1 
\end{array} \right)
\] 
of $\mathrm{SL}_{d+1}(\R)$ whose action on $\R^{d+1}$ is not irreducible.

\medskip
A first important case of reducible actions on $\V$ is when the support of $\mu$ generates a group
\[
\left\{ \left( \begin{array}{cc}
\Gamma_1 & \ast \\ 0 & \Gamma_2
\end{array}\right)\right\}
\]
where $\Gamma_i < \mathrm{SL}(\V_i)$ with $\V_1 \oplus \V_2 = \V$ and $\Gamma_i$ acts irreducibly on $\V_i$. Indeed, in this case, we can study the action on $\V_1$ and on $\V_2$ to get the one on $\V$. This motivates the following

\begin{definition}
Let $\mu$ be a borelian probability measure on $\mathrm{GL}(\V)$.

We say that some subspace $\W$ of $\V$ is adapted to $\mu$ if it is proper, invariant by the subgroup of $\mathrm{GL}(\V)$ spanned by the support of $\mu$ and if there are $\Delta_1>\Delta_2 \in \R$ such that
\begin{itemize}
\item For any $x\in \V/\W \setminus\{0\}$,
\[
\frac 1 n \ln\|g_n \dots g_1 x\| \xrightarrow\, \Delta_1\text{ a.e.}
\]
\item For any $x\in \W \setminus\{0\}$,
\[
\limsup\frac 1 n \ln\|g_n \dots g_1 x\| \leqslant \Delta_2\text{ a.e.}
\]
\end{itemize}
\end{definition}

This definition is only useful since there is always an adapted subspace.

\begin{theorem}[\cite{FK83} or \cite{Hen84}] \label{theorem:FKH}
Let $\mu$ be a borelian probability measure on $\mathrm{GL}(\V)$ having a moment of order $1$.

Then there is some subspace $\W$ of $\V$ that is adapted to $\mu$.
\end{theorem}

This theorem proves that we can always block-triangularize the group $\G_\mu$ spanned by the support of $\mu$. Indeed, we can find by induction $\Delta_1>\dots >\Delta_s \in \R$ and a flag $\V=\W_1 \supset \dots \supset \W_{s+1}:=0$ adapted to $\mu$ : $\W_i$ is $\G_\mu-$invariant and for any $i$ and any $x\in \W_i / \W_{i+1} \setminus\{0\}$,
\begin{equation} \label{equation:drapeau_invariant}
\frac 1n \ln \|g_n \dots g_1 x\| \xrightarrow\, \Delta_i \text{ a.e.}
\end{equation}
The flag that we obtain in this way is included in the flag given by Oseledec's theorem (it means that $\{\W_i\} \subset \{\V_i\}$ and that $\{\Delta_i\} \subset \{\Lambda_i\}$ with $\Delta_1 = \Lambda_1$) but it has the advantage of being invariant and the convergence in equation~\eqref{equation:drapeau_invariant} gives us Furstenberg's law of large numbers without any irreducibility assumption.

\medskip
We are going to precise the convergence in equation~\eqref{equation:drapeau_invariant} through a non-concentration inequality that we state in next
\begin{theorem} \label{theoreme:LDP_produits_reductible}
Let $\mu$ be a borelian probability measure on $\mathrm{GL}(\V)$ having an exponential moment.

Then, for any $\varepsilon \in \R_+^\ast$, there are $C,t\in \R_+^\ast$ such that for any $n\in \N$,
\[
\mu^{\ast n} \left( \left\{ g \in \mathrm{GL}(\V)\middle| \left| \frac 1 n \ln \|g\| - \Lambda_1 \right| \geqslant \varepsilon\right\}\right) \leqslant Ce^{-tn}
\]
Moreover, if $\W$ is adapted to $\mu$ then, for any $x\in \V$ and any $n\in \N$,
\[
\mu^{\ast n} \left(\left\{ g\in \mathrm{GL}(\V)\middle|  e^{(\Lambda_1 - \varepsilon)n} d(x,\W) \leqslant \|gx\| \leqslant e^{(\Lambda_1 + \varepsilon)n} \|x\| \right\}\right) \geqslant 1- Ce^{-tn}
\]
where we noted
\[
d(x,\W) = \inf_{y\in \W} \|x-y\|
\]
\end{theorem}

To prove this theorem, we first prove the following
\begin{lemma} \label{lemma:LDPmatrices}
Let $\mu$ be a borelian probability measure on $\mathrm{GL}(\V)$ having an exponential moment.

Let $\Delta_1>\dots>\Delta_s\in \R$ and $\V:=\W_1 \supset \dots \supset \W_{s+1}:=\{0\}$ be the flag adapted to $\mu$ and given by induction by theorem~\ref{theorem:FKH}.

Then for any $\varepsilon \in \R_+^\ast$, there are $C,t\in \R_+^\ast$ such that for any $n\in \N$ and any $x\in \V\setminus\{0\}$,
\[
\mu^{\ast n} \left( \left\{ g \in \mathrm{GL}(\V)\middle| \Delta_s - \varepsilon \leqslant\frac 1 n \ln \frac{ \|gx\|}{\|x\|}  \leqslant \Delta_1 + \varepsilon\right\}\right) \geqslant 1 - Ce^{-tn}
\]
\end{lemma}

\begin{proof}
The proof of this lemma is an adaptation of the proof of proposition~3.2 in~\cite{BQlin} where Benoist and Quint only have polynomial moments.

First, for $g\in \G$ and $X=\R x\in \prob(\R^d)$, we set
\[
\sigma(g,X) = \ln \frac{ \|gx\|}{\|x\|},\;\varphi(X) = \int_\G \sigma(g,X) \di\mu(g)  \text{ et }\sigma'(g,X) = \sigma(g,X) - \varphi(X)
\]
Then, for any $X=\R x \in \prob(\R^d)$ and any sequence $(g_n) \in \mathrm{GL}(\V)^\N$, noting $X_k = g_k \dots g_1 X$, we have that
\[
\ln \frac{\|g_n \dots g_1 x\|}{\|x\|} = \sum_{k=0}^{n-1} \sigma(g_{k+1}, X_k) = \sum_{k=1}^n \sigma'(g_{k+1}, X_k) + \sum_{k=0}^{n-1} \varphi(X_k)
\]
Now, let $M_n = \sum_{k=0}^{n-1}\sigma'(g_{k+1},X_k)$.

We can compute
\[
\esp\left[ M_{n+1} - M_n \middle| X_0, \dots,X_n \right] = \esp \left[\sigma'(g_{n+1}, X_n) \middle| X_0, \dots ,X_n \right] = \int_\G \sigma'(g,X_n) \di\mu(g) = 0
\]
This proves that $M_n$ is a martingale.
Moreover, for any $\varepsilon \in \R_+^\ast$,
\[
\esp \left[ e^{\varepsilon \left|M_{n+1} - M_n\right|} \middle| X_0, \dots, X_n \right] = \int_\G e^{\varepsilon \left|\sigma'(g,X_n)\right|} \di\mu(g) 
\]
And the inequality
\[
-\ln \|g^{-1} \| \leqslant \sigma(g,X) = \ln \frac{ \|gx\|}{\|x\|} \leqslant \ln \|g\|
\]
shows that
\[
\esp \left[ e^{\varepsilon \left|M_{n+1} - M_n\right|} \middle| X_0, \dots, X_n \right] \leqslant \int_\G e^{2\varepsilon \ln \max(\|g\|, \|g^{-1}\|)} \di\mu(g)
\]
And so, as $\mu$ has an exponential moment, there is $\varepsilon \in \R_+^\ast$ and some constant $C_0$ such that for any $X\in \prob(\V)$ and any $n\in \N$,
\[
\esp_X \left[ e^{\varepsilon \left|M_{n+1} - M_n\right|} \middle| X_0, \dots, X_n \right] \leqslant C_0 \text{ a.e.}
\]
Thus, according to the non-concentration inequality for martingales (see theorem~1.1 in~\cite{LW09}), for any $\varepsilon \in \R_+^\ast$, there is $C,t\in \R_+^\ast$ such that for any $X\in \prob(\V)$ and any $n\in \N$,
\[
\prob_X( |M_n| \geqslant \varepsilon n)\leqslant Ce^{-tn}
\]
(it is important to remark that the constants $C,t$ don't depend on $X$ but only on $\varepsilon$ and $C_0$).

To conclude, we only need to study $\sum_{k=0}^{n-1} \varphi(X_k)$. To do so, we just proved (using Borel-Cantelli-s theorem) that
\[
\frac 1 n M_n \xrightarrow\, 0 \text{ a.e.}
\]
Moreover, by definition of $\Delta_1, \dots , \Delta_s$ and $\W_1, \dots ,\W_s$, we have that for any $X\in \prob(\V)$,
\[
\Delta_s \leqslant \liminf_n \frac 1 n \sum_{k=0}^{n-1} \sigma(g_{k+1}, X_k) \leqslant \limsup_n \frac 1 n \sum_{k=0}^{n-1} \sigma(g_{k+1}, X_k) \leqslant \Delta_1 \text{ a.e.}
\]
This proves that
\[
\Delta_s \leqslant \liminf_n \frac 1 n \sum_{k=0}^{n-1} \varphi(X_k) \leqslant \limsup_n \frac 1 n \sum_{k=0}^{n-1} \varphi(X_k) \leqslant\Delta_1 \text{ a.e.}
\]
And so, for any stationary probability measure $\nu$ on $\prob(\V)$,
\[
\Delta_s \leqslant\int \varphi \di\nu = \lim_n \int \esp_x \frac 1 n \sum_{k=0}^{n-1} \varphi(X_k) \di\nu(x) \leqslant \Delta_1
\]
Finally, using proposition~3.1 of~\cite{BQlin}, we get that for any $\varepsilon \in \R_+^\ast$, there are $C,t\in \R_+^\ast$ such that for any $n\in \N$ and any $X\in \prob(\R^d)$,
\[
\prob_X\left( \Delta_s - \varepsilon \leqslant\frac 1 n\sum_{k=0}^{n-1} \varphi(X_k) \leqslant \Delta_1+ \varepsilon  \right)\leqslant Ce^{-tn}
\]
And this finishes the proof of the lemma.
\end{proof}

\begin{proof}[End of the proof of theorem~\ref{theoreme:LDP_produits_reductible}]
We refer to the proof of proposition~4.1 in~\cite{BQlin}.

First, remark that for any basis $(v_1 , \dots , v_d)$ of $\V$ there is a constant $C$ such that for any $g\in \mathrm{GL}(\V)$,
\[
\frac 1 C\|gv_1\| \leqslant \|g\| \leqslant C \max_{i\in [1,d]} \|gv_i\|
\]
And so, the non-concentration inequality for $\ln\frac{\|gv\|}{\|v\|}$ for any $v\in \V\setminus\{0\}$ implies the one for $\ln\|g\|$.

Then, remark that according to lemma~\ref{lemma:LDPmatrices}, for any $\varepsilon \in \R_+^\ast$, there are $C,t\in \R_+^\ast$ such that for any $x\in \V$ and any $n\in \N$,
\[
\mu^{\ast n} \left(\left\{ g\in \mathrm{GL}(\V) \middle| \|gx\| > e^{(\Lambda_1 + \varepsilon)n} \|x\|\right\}\right) \leqslant Ce^{-tn}
\]

We endow $\V/\W$ with the norm
\[
\|x\|_{\V/\W} := \inf_{y\in \W} \|x-y\|
\]
Then, we have that for any $x\in \V$,
\[
\|x\| \geqslant d(x,\W) = \|\pi x\|_{\V/\W}
\]
where $\pi$ is the projection onto $\V/\W$ and, by definition, since $\W$ is adapted to $\mu$, for any $x\in \V/\W \setminus\{0\}$,
\[
\frac 1 n \ln \|g_n \dots g_1 x\|_{\V/\W} \xrightarrow\, \Lambda_1 \text{ a.e.}
\]
This proves that $\{0\}$ is adapted to the image of $\mu$ in $\mathrm{GL}(\V/\W)$. And so, according to lemma~\ref{lemma:LDPmatrices}, for any $\varepsilon \in \R_+^\ast$, there are $C,t\in \R_+^\ast$ such that for any $n\in \N$ and any $x\in \V/\W \setminus\{0\}$,
\[
\mu^{\ast n }\left(\left\{ g\middle| \left| \frac 1 n \ln \frac{\|gx\|}{\|x\|} - \Lambda_1\right| \geqslant \varepsilon\right\}\right) \leqslant Ce^{-tn}
\]
And this proves that for any $n\in \N$ and any $x\in \V$,
\[
\mu^{\ast n} \left(\left\{ g\middle| \|gx\| \leqslant e^{(\Lambda_1 - \varepsilon)n} d(x,\W) \right\}\right) \leqslant Ce^{-tn}
\]
And this finishes the proof of the theorem.
\end{proof}

We end this section with the study, for a borelian probability measure $\mu$ on $\mathrm{SL}_d(\R) \ltimes \R^d$, of the translation part $b$ of the $\mu^{\ast n}-$generic elements $g=(a,b)$. The aim is to prove that if we make some assumptions on $\mu$, for any $n$, with large $\mu^{\ast n}-$probability, an element $g=(a,b)\in \mathrm{SL}_d(\R) \ltimes \R^d$ is such that $\|b\|\approx e^{\Lambda_1 n}$.

To do so, we wirst compare the Lyapunov exponents of $\mu$ and of it's projection on $\mathrm{SL}_d(\R)$ in next

\begin{lemma}\label{lemma:lyapunov_affine}
Let $\mu$ be a borelian probability measure on $\mathrm{SL}_d(\R) \ltimes \R^d$ having an exponential moment and let $\mu_0$ be the projection of $\mu$ onto $\mathrm{SL}_d(\R)$.

See $\mu$ as a probability measure on $\mathrm{SL}_{d+1}(\R)$ and define $\Lambda_1(\mu)$ this way.

Then,
\[
\Lambda_1(\mu) = \Lambda_1(\mu_0)
\]
\end{lemma}

\begin{proof}
First for any $n\in \N^\ast$ and any $\varepsilon \in \R_+^\ast$,
\[
\mu^{\otimes \N} \left(\left\{  \left| \frac 1 n \ln \|g_{2n} \dots g_{n+1}\| - \Lambda_1 (\mu)\right| \geqslant \varepsilon \right\}\right) = \mu^{\otimes \N} \left(\left\{  \left|\frac 1 n \ln \|g_n \dots g_1\| - \Lambda_1 (\mu)\right| \geqslant \varepsilon \right\}\right)
\]
So, according to lemma~\ref{theoreme:LDP_produits_reductible},
\[
\frac 1 n \ln \| g_{2n} \dots g_{n+1} \| \xrightarrow\, \Lambda_1(\mu) \;\mu^{\otimes\N}-\text{a.e.}
\]

We can write, for any $g\in \mathrm{SL}_d(\R) \ltimes \R^d$,
\[
g =\left( \begin{array}{cc}a & b\\0&1 \end{array}\right)
\]
and if $g$ has law $\mu$, then $a$ has law $\mu_0$.

Thus,
\[
g_{2n} \dots g_1 = \left(\begin{array}{cc}
a_{2n} \dots a_1 &a_{2n} \dots a_{n+1} \sum_{k=1}^{n} a_n \dots a_{k+1} b_k + \sum_{k=n+1}^{2n} a_{2n} \dots a_{k+1} b_k \\ 0 & 1
\end{array} \right)
\]
And so,,
\[
\|g_n \dots g_1 \| \geqslant \|a_n \dots a_1\|
\]
This proves that $\Lambda_1(\mu) \geqslant \Lambda_1(\mu_0)$.

Let now $\Omega \in \G^\N$ be such that $\mu^{\otimes \N} (\Omega)=1$ and for any $(g_n) \in \Omega$,
\[
\lim_n \frac 1 n \ln \|g_{n} \dots g_1 \| =\lim_n \frac 1 n \ln \|g_{2n} \dots g_{n+1}\| = \Lambda_1(\mu)
\]
and
\[
\lim_n\frac 1 n \ln \|a_{n} \dots a_1 \| = \lim_n \frac 1n \ln \|a_{2n} \dots a_{n+1}\| = \Lambda_1(\mu_0)
\]
This way, for any $\varepsilon \in \R_+^\ast$, there is $N_\varepsilon$ such that for any $n\in \N$ with $n\geqslant N_\varepsilon$,
\[
e^{(\Lambda_1(\mu) - \varepsilon)n} \leqslant\|g_n \dots g_1 \| , \|g_{2n} \dots g_{n+1} \| \leqslant e^{(\Lambda_1(\mu) + \varepsilon)n}
\]
and,
\[\|a_n \dots a_1 \|, \|a_{2n} \dots a_{n+1} \| \leqslant e^{(\Lambda_1(\mu_0) + \varepsilon)n}
\]
Thus, we have that for any $\varepsilon\in \R_+^\ast$ and any large enough $n$,
\[
e^{2 (\Lambda_1(\mu)-\varepsilon)n} \leqslant \|g_{2n} \dots g_1 \| \leqslant \max\left( e^{2(\Lambda_1(\mu_0) + \varepsilon)n}, e^{(\Lambda_1(\mu) + \Lambda_1(\mu_0) + 2\varepsilon)n} + e^{(\Lambda_1(\mu) + \varepsilon)n} \right)
\]
This proves that for any $\varepsilon \in \R_+^\ast$,
\[
2(\Lambda_1(\mu) - \varepsilon) \leqslant \max (2\Lambda_1(\mu_0) + \varepsilon, \Lambda_1(\mu) + \Lambda_1(\mu_0) + 2\varepsilon)
\]
And so, we get that
\[
2\Lambda_1(\mu) \leqslant \max( 2\Lambda_1(\mu_0), \Lambda_1(\mu) + \Lambda_1(\mu_0))
\]
Finally, as we already proved that $\Lambda_1(\mu_0) \leqslant \Lambda_1(\mu)$, the previous inequality actually is an equality and we get the expected result.
\end{proof}

Until now, we didn't say anything on the positivity of $\Lambda_1$.

If $\mu$ is a measure on $\mathrm{SL}(\V)$, then, as for $\mu-$a.e. $g\in \G$, $\det(g)=1$, we have that $\lambda_1 + \dots + \lambda_{\dim(\V)} = 0$ and so, $\lambda_1=0$ if and only if for any $i$, $\lambda_i=0$.

To get conditions that ensure that $\lambda_1>0$, we will say that a subgroup $\Hb$ of $\mathrm{SL}_d(\R)$ is strongly irreducible if it doesn't fix any non trivial finite union of subspaces of $\R^d$.

We remind the following result
\begin{theorem}[see~\cite{BL85} or~\cite{Gu81}] \label{theorem:positivite_Lyapunov}
Let $\mu$ be a borelian probability measure on $\mathrm{SL}(\V)$ having a moment of order $1$ and such that the subgroup spanned by the support of $\mu$ is strongly irreducible and non-compact.
\[
\int_{\mathrm{SL}_d(\R)} \ln \|g\|\di\mu(g) <+\infty
\]

Then, $\Lambda_1 >0$.  
\end{theorem}

\begin{corollary}\label{corollaire:longueurs_translations}
Let $\mu$ be a borelian probability measure on $\mathrm{SL}_d(\R) \ltimes \R^d$ having an exponetial moment and that is not concentrated on $\mathrm{SL}_d(\R) \ltimes\{0\}$. Assume that the projection onto $\mathrm{SL}_d(\R)$ of the subgroup spanned by the support of $\mu$ is strongly irreducible and non-compact.

Then, for any $\varepsilon \in \R_+^\ast$, there are $C,t\in \R_+^\ast$ such that for any $n\in \N$ and any $x\in\R^d$,
\[
\mu^{\ast n} \left(\left\{ g \in \mathrm{SL}_d(\R) \ltimes \R^d \middle|  \left| \frac 1 n \ln\frac{1+\|gx\|}{1+\|x\|} - \Lambda_1 \right| \geqslant \varepsilon \right\}\right) \leqslant Ce^{-tn}
\]
\end{corollary}

\begin{proof}
We see $\mu$ as a probability measure on $\mathrm{SL}_{d+1}(\R)$ and we are going to prove that the subspace $\W$ adapted to $\mu$ and given by theorem~\ref{theorem:FKH} is $\{0\}$.

First, since the projection of $\supp\mu$ on $\mathrm{SL}_d(\R)$ spans a group $\G_\mu$ that acts strongly irreducibly and non-compactly on $\R^d$, the only subspaces of $\R^{d+1}$ that can be invariant by $\supp\mu$ are $\R^{d+1}$, $\mathrm{Vect}(e_1, \dots e_{d})$, $\mathrm{Vect}(e_{d+1})$ and $\{0\}$. But, assuming that $\supp\mu$ is not a subset of $\mathrm{SL}_d(\R) \ltimes\{0\}$ implies that $\mathrm{Vect}(e_{d+1})$ is not invariant by the group spans by the support of $\mu$. Moreover, according to theorem~\ref{theorem:positivite_Lyapunov}, there is $\Lambda_1\in \R_+^\ast$ such that for any $x\in \mathrm{Vect}(e_1,\dots, e_d) \setminus\{0\}$,
\[
\frac 1 n \ln \|g_n \dots g_1 x\| \xrightarrow\, \Lambda_1 \text{ a.e.}
\]
Then, for any $x\in \R^{d+1} \setminus \mathrm{Vect}(e_1, \dots, e_{d+1})$, we have, according to lemma~\ref{lemma:lyapunov_affine} that
\[
\limsup_n \frac 1 n \ln \|g_n \dots g_1 x\| \leqslant \Lambda_1
\]
So we also have that $\W \not=\mathrm{Vect}(e_1, \dots ,e_d)$ and, as $\W$ is proper, we have that $\W=\{0\}$.

Thus, according to lemma~\ref{theoreme:LDP_produits_reductible}, for any $\varepsilon \in\R_+^\ast$, there are $C,t\in \R_+^\ast$ such that for any $x\in \R^{d+1} \setminus\{0\}$,
\[
\mu^{\ast n} \left( \left\{g \in \mathrm{SL}_d(\R) \ltimes \R^d \middle| \left| \frac 1 n \ln \frac{\|gx\|}{\|x\|}- \Lambda_1 \right| \geqslant \varepsilon \right\}\right) \leqslant Ce^{-tn}
\]
And in particular, with $x=x_0+e_{d+1}$ for $x_0\in \mathrm{Vect}(e_1, \dots, e_d)$ we get the expected result.
\end{proof}

\bibliographystyle{amsalpha}
\bibliography{biblio}

\end{document}